 \documentclass[preprint,review,12pt]{elsarticle}
\usepackage{graphicx}

\usepackage{amssymb,amsmath}
\usepackage{amsthm}

\journal{}

\begin{document}

\begin{frontmatter}



\title{Square-mean weighted pseudo almost automorphic solutions for nonautonomous stochastic differential equations driven by L$\acute{e}$vy noise\tnoteref{t1}}

\author[]{Kexue Li}
\ead{kexueli@gmail.com}

\tnotetext[t1]{This work is partially supported by the National Natural Science Foundation of China under contacts No. 11201366 and 11131006.}

\address{School of Mathematics and Statistics, Xi'an Jiaotong University, Xi'an 710049, PR China}

\begin{abstract}In this paper, we introduce the concepts of Poisson square-mean almost automorphy and Poisson square-mean weighted pseudo almost automorphy. Using the theory of evolution family and stochastic analysis techniques, we establish the existence and uniqueness results of square-mean weighted pseudo almost automorphic solutions for some linear and semilinear  nonautonomous stochastic differential equations driven by L$\acute{e}$vy noise. Moreover, we investigate the global and local exponentially stability of the square-mean weighted pseudo almost automorphic solutions.

\end{abstract}


\end{frontmatter}
\section{Introduction}
Almost periodic functions and almost periodic solutions of deterministic differential equations had been extensively investigated \cite{AF}, \cite{BV}. In the real world, stochastic perturbation is unavoidable. Almost periodic solutions of stochastic differential equations were studied by several authors, Tudor \cite{T} and Da prato \cite{DC} studied the almost periodic solution in the distribution sense for stochastic evolution equations. Bezandry \cite{Bezandry} studied the existence of square-mean almost periodic solutions to some semilinear stochastic equations. For more details about almost periodic solutions of stochastic differential equations, we refer to \cite{BD}.

The concept of almost automorphy is introduced by Bochner \cite{Bochner} and is related to some aspects of differential geometry. The almost automorphy is a generalization of almost periodicity.  Almost automorphic solutions to deterministic ordinary as well as abstract differential equations have been intensively studied \cite{Goldstein, Ezzinbi, Diagana, CE, Cara, Caraballo}.  In the book \cite{GMN}, the author gave an overview about the theory of almost automorphic functions and their applications to deterministic differential equations. Fu and Liu \cite{FM} introduced the concept of square-mean almost automorphy for stochastic processes and they studied square-mean almost automorphic solutions to some linear and non-linear stochastic differential equations. Cao et al \cite{Cao} introduced the concept of $p$-mean almost automorphy for stochastic processes. The existence, uniqueness and exponential stability of quadratic-mean almost automorphic mild solutions for a class of stochastic functional differential equations were discussed.

Zhang \cite{Zhang} introduced the pseudo almost periodic function, which is a generalization of the almost periodic function. Some authors studied pseudo-almost periodic solutions of deterministic differential equations, see \cite{Yuan,EH,Diag}. We refer to \cite{CM} for the pseudo almost periodic solutions of stochastic differential equations. N'Gu$\acute{e}$r$\acute{e}$kata \cite{GMN} introduced the pseudo almost automorphic function, which is a generalization of the pseudo almost periodic function. Xiao et al \cite{Xiao} obtained an existence and uniqueness theorem of pseudo almost automorphic mild solutions to semilinear differential equations in Banach spaces. Blot et al \cite{Blot} proposed the notion of weighted pseudo almost automorphic functions, which generalizes the weighted pseudo almost periodicity \cite{TD}, and the authors established the existence and uniqueness theorem for pseudo almost automorphic mild solutions to semilinear differential equations in a Banach space.  Liu and Song \cite{Liu} proved the existence and uniqueness of a weighted pseudo almost automorphic mild solution of the semilinear nonautonomous evolution equation. Chen and Lin \cite{Chen} introduced the concept of square-mean pseudo almost automorphy, and the authors investigated the existence and uniqueness and the global stability of the square-mean pseudo almost automorphic solutions for a class of stochastic evolution equations.

It should be noted that many stochastic processes may not be square-mean pseudo almost automorphic since time average terms do not converge to zero as the time goes to positive infinity. Chen and Lin \cite{CZ} introduced the concept of the square-mean weighted pseudo almost automorphy, which is a generalization of the square-mean pseudo almost automorphy.  The authors studied the well-posedness of the square-mean weighted pseudo almost automorphic solutions for a general class of non-autonomous stochastic evolution equations that satisfy either global or local Lipschitz condition.  Moreover, they estimated the boundedness of attractive domain for the case where the local Lipschitz condition is satisfied.

So far, most of the studies on almost automorphic solutions for stochastic differential equations are concerned with equations perturbed by Brownian motion. However, many real models involve jump perturbations, or more general L$\acute{e}$vy noise. Wang and Liu \cite{YZ} first introduce the concept of Poisson square-mean almost periodicity and study the existence, uniqueness and stability of square-mean almost period solutions for stochastic evolution equations driven by L$\acute{e}$vy noise. In this paper, we introduce the concepts of Poisson square-mean almost automorphy and Poisson square-mean weighted pseudo almost automorphy. We study the Poisson square-mean weighted pseudo almost automorphic solutions for nonautonomous stochastic differential equations driven by L$\acute{e}$vy noise.

This paper is organized as follows. In Section 2, we recall the L$\acute{e}$vy process and introduce the concepts of some abstract spaces, the
Poisson square-mean almost automorphy and Poisson square-mean weighted pseudo almost automorphy and study some of their properties. In Section 3, we investigate square-mean almost automorphic solutions of nonautonomous linear stochastic differential equations under some suitable conditions.
Based on the results in Sections 2 and 3, we prove the existence and uniqueness of square-mean weighted pseudo almost automorphic solutions of nonautonomous linear stochastic differential equations is obtained in Section 4. In Section 5, we prove the existence and uniqueness of square-mean weighted pseudo almost automorphic solutions of nonautonomous semilinear stochastic differential equations. In Section 6, we discuss the exponential stability of the unique square-mean almost automorphic solution in square-mean sense. Finally, we consider the square-mean weighted pseudo almost automorphic solution for a stochastic heat equation.

\section{Preliminaries}
Throughout this paper, $(H,\|\cdot\|)$ and $(V,|\cdot|)$ are assumed to be real separable Hilbert spaces. Let $(\Omega, \mathcal{F}, P)$ be a probability space with given filtration $(\mathcal{F}_{t})_{t\geq 0}$ satisfying the usual condition. $L(V,H)$ denotes the space of all bounded linear operators from $V$ to $H$. $\mathcal{L}^{2}(P,H)$ stands for the space of all $H$-valued  random variable $Y$ such that $\mathbb{E}\|Y\|^{2}=\int_{\Omega}\|Y\|^{2}dP< \infty$, which is a Banach space equipped with the norm $\|Y\|_{2}:=(\mathbb{E}\|Y\|^{2})^{\frac{1}{2}}$. We  consider a L$\acute{e}$vy process with values in $V$.
\subsection{L$\acute{e}$vy process}
In this subsection, we recall the definition and some properties of L$\acute{e}$vy processes. \\
 \textbf{Definition 2.1.} (see \cite{DA, PZ}.) A $V$-valued stochastic process $L=(L(t), t\geq 0)$ is called a L$\acute{e}$vy process if\\
(i) $L(0)=0$ almost surely;\\
(ii)  $L$ has independent and stationary increments; \\
(iii) $L$ is stochastic continuous, i.e., for all $\epsilon>0$ and $s>0$,
\begin{align*}
\lim_{t\rightarrow s}P(|L(t)-L(s)|_{V}>\epsilon)=0;
\end{align*}
Every L$\acute{e}$vy process is c$\grave{a}$dl$\grave{a}$g. Given a L$\acute{e}$vy process $L$, we define the process of jumps of $L$ by
$\Delta L(t)=L(t)-L(t-), t\geq 0$, where $L(t-)=\lim_{s\uparrow t}L(s)$.  For any Borel set $B$ in $V-\{0\}$, define the counting measure
\begin{align*}
N(t,B)(\omega)=\sharp\{0\leq s\leq t: \Delta L(s)(\omega)\in B\}=\sum_{0\leq s\leq t} \chi_{B}(\triangle X(s))(\omega)),
\end{align*}
where $\sharp$ means the counting and $\chi_{B}$ is the indicator function. We write $\nu(\cdot)=E(N(1,\cdot))$  and call it the intensity measure associated with $L$.
If a Borel set $B$ in $V-\{0\}$ is bounded below (that is, $0\notin \overline{B}$, where $\overline{B}$ is the closure of $B$), then
$N(t,B)< \infty$ almost surely for all $t\geq 0$ and $\{N(t,B), t\geq 0\}$ is a Poisson process with the intensity $\nu(B)$. $N$ is called the Poisson random measure. For each $t\geq 0$ and $B$ bounded below, the compensated Poisson random measure is defined by
\begin{align*}
\widetilde{N}(t,B)=N(t,B)-t\nu(B).
\end{align*}
 \textbf{Proposition 2.2.} (L$\acute{e}$vy-It$\hat{o}$ decomposition) If $L$ is a $V$-valued L$\acute{e}$vy process, then there exist $a\in U$, a $V$-valued Wiener process $W$ with convariance operator $Q$, and an independent Poisson random measure on $R^{+}\times (V-\{0\})$ such that for each $t\geq 0$,
 \begin{align}\label{levy}
 L(t)=at+W(t)+\int_{|x|_{V}<1}x\widetilde{N}(t,dx)+\int_{|x|_{V}\geq 1}x N(t,dx),
 \end{align}
 where the Poisson random measure $N$ has the intensity measure $\nu$ which satisfies
 \begin{align}\label{measure}
 \int_{V}(|y|_{V}^{2}\wedge1)\nu(dy)<\infty,
  \end{align}
and $\widetilde{N}$ is the compensated Poisson random measure of $N$.

Given two independent, identically distributed L$\acute{e}$vy processes $L_{1}$ and $L_{2}$ with decompositions as in Proposition 2.2, let
\begin{equation*}
 L(t)=\ \left\{\begin{aligned} &L_{1}(t), \ t\geq 0,\\
& L_{2}(-t)\ t\leq 0.
\end{aligned}\right.
\end{equation*}
Then $L$ is a two-sided L$\acute{e}$vy process. In this paper, we need consider two-sided L$\acute{e}$vy processes. The two-sided L$\acute{e}$vy process
$L$ is defined on the filtered probability space  $(\Omega, \mathcal{F}, P, (\mathcal{F}_{t})_{t\in R})$. We assume that the covariance operator $Q$ of $W$ is of trace class, i.e., $\mbox{Tr} Q<\infty$.\\
\textbf{Remark 2.3}(see \cite{YZ}.) It follows from (\ref{measure}) that $\int_{|x|_{V}\geq 1}\nu(dx)<\infty$. For convenience, we denote
\begin{align*}
c: =\int_{|x|_{V}\geq 1}\nu(dx)
\end{align*}
throughout the paper.
\subsection{Poisson square-mean weighted pseudo almost automorphic process}
In this subsection, we present some preliminaries for later use. \\
 \textbf{Definition 2.4.} A stochastic process $Y: \mathbb{R}\rightarrow \mathcal{L}^{2}(P,H)$ is said to be $\mathcal{L}^{2}$-bounded if there exists a constant $M>0$ such that
  \begin{equation*}
  \mathbb{E}\|Y(t)\|^{2}=\int_{\Omega} \|Y(t)\|^{2}dP\leq M.
 \end{equation*}
 \textbf{Definition 2.5.} A stochastic process $Y: \mathbb{R}\rightarrow \mathcal{L}^{2}(P,H)$ is said to be $\mathcal{L}^{2}$-continuous if for any $s\in \mathbb{R}$,
 \begin{equation*}
 \lim_{t\rightarrow s}\mathbb{E}\|Y(t)-Y(s)\|^{2}=0.
 \end{equation*}

 Denote by $SBC(\mathbb{R}, \mathcal{L}^{2}(P,H))$ the collection of all the $\mathcal{L}^{2}$-bounded and $\mathcal{L}^{2}$-continuous processes. \\
 \textbf{Remark 2.6.} $SBC(\mathbb{R}, \mathcal{L}^{2}(P,H))$ is a Banach space equipped with the norm $\|Y\|_{\infty}=\sup _{t\in \mathbb{R}}(\mathbb{E}\|Y(t)\|^{2})^{\frac{1}{2}}$. \\

   Let $\mathcal{U}$ be the set of all functions which are positive and locally integrable over $\mathbb{R}$. For given $r>0$ and $\rho\in \mathcal{U}$, define $m(r,\rho)=\int_{-r}^{r}\rho(t)dt$ and
 \begin{equation*}
  \mathcal{U}_{\infty}=\{\rho\in \mathcal{U}|\lim_{r\rightarrow +\infty}m(r,\rho)=+\infty\}.
  \end{equation*}
  \textbf{Definition 2.7.} (see \cite{CZ}.) For $\rho\in \mathcal{U}_{\infty}$, define a class of stochastic processes
  \begin{equation*}
  SBC_{0}(\mathbb{R},\rho)=\{Y\in SBC(\mathbb{R}, \mathcal{L}^{2}(P,H))|\lim_{r\rightarrow +\infty}\frac{1}{m(r,\rho)}\int_{-r}^{r}E\|Y(t)\|^{2}\rho(t)dt=0\}.
  \end{equation*}
  \textbf{Remark 2.8.} (see \cite{CZ}.) $SBC_{0}(\mathbb{R},\rho)$ is a linear closed subspace of $SBC(\mathbb{R}, \mathcal{L}^{2}(P,H))$. \\
   \textbf{Remark 2.9.} (see \cite{CZ}.) $SBC_{0}(\mathbb{R},\rho)$ equipped with the norm $\|Y\|_{\infty}$ is a Banach space. \\
  \textbf{Definition 2.10.}(see \cite{FM}.) An $\mathcal{L}^{2}$-continuous stochastic process $Y: \mathbb{R}\rightarrow \mathcal{L}^{2}(P,H)$ is said to be square-mean almost automorphic if every sequence of real numbers $\{s'_{n}\}$ has a subsequence $\{s_{n}\}$ such that for some stochastic process $\widetilde{Y}:\mathbb{R}\rightarrow \mathcal{L}^{2}(P,H)$, $\lim_{n\rightarrow\infty}\mathbb{E}\|Y(t+s_{n})-\widetilde{Y}(t)\|^{2}=0$ and  $\lim_{n\rightarrow\infty}\mathbb{E}\|\widetilde{Y}(t-s_{n})-Y(t)\|^{2}=0$ hold for each $t\in \mathbb{R}$.

  The collection of all square-mean almost automorphic process $Y: \mathbb{R}\rightarrow \mathcal{L}^{2}(P,H)$ is denoted by $SAA(\mathbb{R}, \mathcal{L}^{2}(P,H))$. It is a Banach space with the norm $\|Y\|_{\infty}$. \\
  \textbf{Remark 2.11.} (see \cite{FM}.) If $Y\in SAA(\mathbb{R}, \mathcal{L}^{2}(P,H))$, then $Y$ is bounded, that is, $\|Y\|_{\infty}<\infty$. \\
\textbf{Proposition 2.12.} (see \cite{FM}.) Let $f: \mathbb{R}\times \mathcal{L}^{2}(P,H)\rightarrow \mathcal{L}^{2}(P,H)$, $(t,Y)\mapsto f(t,Y)$ be square-mean almost automorphic in $t\in \mathbb{R}$ for each $Y\in \mathcal{L}^{2}(P,H)$, and assume that $f$ satisfies the Lipschitz condition in the following sense:
    \begin{equation*}
    \mathbb{E}\|f(t,Y)-f(t,Z)\|^{2}\leq L\|Y-Z\|^{2}
    \end{equation*}
    for all $Y,Z\in \mathcal{L}^{2}(P,H)$, and for each $t\in \mathbb{R}$, where $L>0$ is independent of $t$. Then for any almost automorphic process $Y: \mathbb{R}\rightarrow \mathcal{L}^{2}(P,H)$, the stochastic process $F:\mathbb{R}\rightarrow \mathcal{L}^{2}(P,H)$ given by $F(t):=f(t,Y(t))$ is square-mean almost automorphic. \\
\textbf{Definition 2.13.} (see \cite{CZ}.)  An $\mathcal{L}^{2}$-continuous stochastic process $f: \mathbb{R}\rightarrow\mathcal{L}^{2}(P,H)$ is said to be square-mean weighted pseudo almost automorphic with respect to $\rho\in \mathcal{U}_{\infty}$ if it can be decomposed as $f=g+\varphi$, where $g\in SAA(\mathbb{R}, \mathcal{L}^{2}(P,H))$ and $\varphi\in SBC_{0}(\mathbb{R}, \rho)$.

  The collection of all square-mean weighted pseudo almost automorphic processes with respect to $\rho$ is denoted by $SWPAA(\mathbb{R}, \rho)$. \\
 \textbf{Definition 2.14.} A set $\mathcal{D}$  is said to be translation invariant if for any $f(t)\in \mathcal{D}$, $f(t+\tau)\in \mathcal{D}$ for any $\tau\in \mathbb{R}$.

   Denote $\mathcal{U}^{inv}=\{\rho\in \mathcal{U}_{\infty}|SBC_{0}(\mathbb{R},\rho) \ \mbox{is translation invariant}$\}. \\
   \textbf{Lemma 2.15.} (see \cite{CZ}.) For $\rho\in\mathcal{U}^{inv}$, $SWPAA(\mathbb{R}, \rho)$ equipped with the norm $\|Y\|_{\infty}$ is a Banach space. \\
  Denote
  \begin{equation*}
  SAA(\mathbb{R}\times \mathcal{L}^{2}(P,H), \mathcal{L}^{2}(P,H))=\{g(t,Y)\in SAA(\mathbb{R},\mathcal{L}^{2}(P,H))|\ \mbox{for any}\ Y\in \mathcal{L}^{2}(P,H)\},
  \end{equation*}
and
  \begin{equation*}
  SBC_{0}(\mathbb{R}\times \mathcal{L}^{2}(P,H), \rho)=\{\varphi(t,Y)\in SBC_{0}(\mathbb{R},\rho)|\ \mbox{for any}\ Y\in \mathcal{L}^{2}(P,H)\}.
  \end{equation*}
  \textbf{Definition 2.16.} (see \cite{CZ}.)  An $\mathcal{L}^{2}$-continuous stochastic process $f: \mathbb{R}\times \mathcal{L}^{2}(P,H)\rightarrow\mathcal{L}^{2}(P,H)$ is said to be square-mean weighted pseudo almost automorphic with respect to $\rho\in \mathcal{U}_{\infty}$ in $t$ for any $Y\in \mathcal{L}^{2}(P,H)$ if it can be decomposed as $f=g+\varphi$, where $g\in SAA(\mathbb{R}\times \mathcal{L}^{2}(P,H), \mathcal{L}^{2}(P,H))$ and $\varphi\in SBC_{0}(\mathbb{R}\times \mathcal{L}^{2}(P,H), \rho)$. We denote all such stochastic process by $SWPAA(\mathbb{R}\times \mathcal{L}^{2}(P,H), \rho))$.\\
  \textbf{Lemma 2.17.} (see \cite{CZ}.) Suppose $\rho\in\mathcal{U}^{inv}$, $f(t,x)\in SWPAA(\mathbb{R}\times \mathcal{L}^{2}(P,H), \rho)$, and there exists a constant $L>0$ such
that for any $x,y\in \mathcal{L}^{2}(P,H)$,
\begin{align*}
\mathbb{E}\|f(t,x)-f(t,y)\|^{2}\leq L\mathbb{E}\|x-y\|^{2}.
\end{align*}
Then, for any $x\in SWPAA(\mathbb{R}, \rho)$, we have $f(t,x)\in SWPAA(\mathbb{R}, \rho)$.\\
\textbf{Definition 2.18.} A continuous stochastic process $f: \mathbb{R}\times \mathcal{L}^{2}(P,H)\rightarrow L(V, \mathcal{L}^{2}(P,H))$ is said to be square-mean almost automorphic in $t\in \mathbb{R}$ for all $Y\in \mathcal{L}^{2}(P,H)$ if for every sequence of real numbers $\{s'_{n}\}$, there exists a subsequence $\{s_{n}\}$ such that for some function $\widetilde{f}: \mathbb{R}\times \mathcal{L}^{2}(P,H)\rightarrow L(V, \mathcal{L}^{2}(P,H))$
   \begin{equation*}
  \lim_{n\rightarrow \infty}\mathbb{E}\|f(t+s_{n},Y)-\widetilde{f}(t,Y)\|^{2}_{L(V,\ \mathcal{L}^{2}(P,H))}=0,
  \end{equation*}
  and
  \begin{equation*}
  \lim_{n\rightarrow \infty}\mathbb{E}\|\widetilde{f}(t-s_{n},Y)-f(t,Y)\|^{2}_{L(V,\ \mathcal{L}^{2}(P,H))}=0
  \end{equation*}
for all $Y\in \mathcal{L}^{2}(P,H)$ and each $t\in \mathbb{R}$.

The collection of all square-mean almost automorphic function $f: \mathbb{R}\times \mathcal{L}^{2}(P,H)\rightarrow L(V, \mathcal{L}^{2}(P,H))$ is denoted by $SAA(\mathbb{R}\times \mathcal{L}^{2}(P,H), L(V,\mathcal{L}^{2}(P,H)))$. \\ For $f:\mathbb{R}\times \mathcal{L}^{2}(P,H)\rightarrow L(V, \mathcal{L}^{2}(P,H))$ satisfying the corresponding property in Definition 2.16, we also denote all these processes as
  $SWPAA(\mathbb{R}, \rho)$.\\
\textbf{Definition 2.19.} A stochastic process $J(t,x): \mathbb{R}\times V\rightarrow \mathcal{L}^{2}(P,H)$ is said to be Poisson stochastically bounded if there exists a constant $M>0$ such that
\begin{align*}
\int_{V}\mathbb{E}\|J(t,x)\|^{2}\nu(dx)\leq M,
\end{align*}
for all $t\in \mathbb{R}$.

By $PSB(\mathbb{R}\times V, \mathcal{L}^{2}(P,H))$, we denote the collection of all Poisson stochastically bounded processes.\\
\textbf{Definition 2.20.} A stochastic process $J(t,x): \mathbb{R}\times V\rightarrow \mathcal{L}^{2}(P,H)$ is said to be Poisson stochastically continuous if
\begin{align*}
\lim_{t\rightarrow s}\int_{V}\mathbb{E}\|J(t,x)-J(s,x)\|^{2}\nu(dx)=0.
\end{align*}

By $PSC(\mathbb{R}\times V, \mathcal{L}^{2}(P,H))$, we denote the collection of all Poisson stochastically continuous processes. By $PSBC(\mathbb{R}\times V, \mathcal{L}^{2}(P,H)$, we denote the collection of all Poisson stochastically  bounded and continuous processes. \\
\textbf{Definition 2.21.} Let $F: \mathbb{R}\times  \mathcal{L}^{2}(P,H)\times V\rightarrow \mathcal{L}^{2}(P,H)$, $(t,Y,x)\mapsto F(t,Y,x)$, $F$ is said to be Poisson stochastically continuous, if
\begin{align*}
\int_{V}\mathbb{E}\|F(t,Y,x)-F(t',Y',x)\|^{2}\nu(dx)\rightarrow0 \ \mbox{as} \ (t,Y)\rightarrow (t',Y').
\end{align*}
By $PSC(\mathbb{R}\times \mathcal{L}^{2}(P,H)\times V\rightarrow \mathcal{L}^{2}(P,H))$, resp. $PSBC(\mathbb{R}\times \mathcal{L}^{2}(P,H)\times V\rightarrow \mathcal{L}^{2}(P,H)$, we denote the collection of all Poisson stochastically  continuous processes, resp. Poisson stochastically  bounded and continuous processes. \\
\textbf{Definition 2.22.} A stochastic process $J(t,x)\in PSBC_{0}(\mathbb{R}\times V, \mathcal{L}^{2}(P,H))$, provided that $J(t,x)\in PSBC(\mathbb{R}\times V, \mathcal{L}^{2}(P,H))$ and
\begin{align*}
\lim_{r\rightarrow +\infty}\frac{1}{2r}\int_{-r}^{r}\int_{V}\mathbb{E}\|J(t,x)\|^{2}\nu(dx)dt=0.
\end{align*}
\textbf{Definition 2.23.} For $\rho\in \mathcal{U}_{\infty}$, define a class of stochastic processes
\begin{align*}
  PSBC_{0}(\mathbb{R}\times V, \rho)=\{&J(t,x)\in PSBC(R\times V, \mathcal{L}^{2}(P,H))| \\ &\lim_{r\rightarrow +\infty}\frac{1}{m(r,\rho)}\int_{-r}^{r}\int_{V}\mathbb{E}\|J(t,x)\|^{2}\rho(t)\nu(dx)dt=0\}.
  \end{align*}

   Denote $\mathcal{U}_{p}^{inv}=\{\rho\in \mathcal{U}_{\infty}|PSBC_{0}(\mathbb{R}\times V,\rho) \ \mbox{is translation invariant}$\}. \\
\textbf{Remark 2.24.}  If $\rho\equiv 1$, it is obvious that $PSBC_{0}(\mathbb{R}\times V, \rho)$ reduces to $PSBC_{0}(\mathbb{R}\times V, \mathcal{L}^{2}(P,H))$.\\
\textbf{Definition 2.25.} A stochastic process $J: \mathbb{R}\times V\rightarrow \mathcal{L}^{2}(P,H)$, $(t,x)\mapsto J(t,x)$ is said to be Poisson
  square-mean almost automorphic in $t\in \mathbb{R}$ if $J$ is Poisson stochastically continuous and for every sequence of real numbers $\{s'_{n}\}$, there exists a subsequence $\{s_{n}\}$ such that for some function $\widetilde{J}: \mathbb{R}\times V\rightarrow \mathcal{L}^{2}(P,H)$, $(t,x)\mapsto \widetilde{J}(t,x)$ is Poisson stochastically continuous  such that
  \begin{equation*}
  \lim_{n\rightarrow \infty}\int_{V}\mathbb{E}\|J(t+s_{n},x)-\widetilde{J}(t,x)\|^{2}\nu(dx)=0,
  \end{equation*}
  and
  \begin{equation*}
  \lim_{n\rightarrow \infty}\int_{V}\mathbb{E}\|\widetilde{J}(t-s_{n},x)-J(t,x)\|^{2}\nu(dx)=0
  \end{equation*}
for each $t\in \mathbb{R}$.

The collection of all Poisson square-mean almost automorphic functions $J: \mathbb{R}\times V\rightarrow \mathcal{L}^{2}(P,H)$ is denoted by $PSAA(\mathbb{R}\times V, \mathcal{L}^{2}(P,H))$. \\
\textbf{Definition 2.26.} A stochastic process $F: \mathbb{R}\times \mathcal{L}^{2}(P,H)\times V\rightarrow \mathcal{L}^{2}(P,H)$, $(t,Y,x)\mapsto F(t,Y,x)$ is said to be Poisson square-mean almost automorphic in $t\in \mathbb{R}$ for each $Y\in \mathcal{L}^{2}(P,H)$ if $F$ is Poisson stochastically continuous and for every sequence of real numbers $\{s'_{n}\}$, there exists a subsequence $\{s_{n}\}$ such that for some function $\widetilde{F}: \mathbb{R}\times \mathcal{L}^{2}(P,H)\times V\rightarrow \mathcal{L}^{2}(P,H)$, $(t,Y,x)\mapsto \widetilde{F}(t,Y,x)$ is Poisson stochastically continuous such that
  \begin{equation*}
  \lim_{n\rightarrow \infty}\int_{V}\mathbb{E}\|F(t+s_{n},Y,x)-\widetilde{F}(t,Y,x)\|^{2}\nu(dx)=0,
  \end{equation*}
  and
  \begin{equation*}
  \lim_{n\rightarrow \infty}\int_{V}\mathbb{E}\|\widetilde{F}(t-s_{n},Y,x)-F(t,Y,x)\|^{2}\nu(dx)=0
  \end{equation*}
  for all $Y\in \mathcal{L}^{2}(P,H)$ and each $t\in \mathbb{R}$. The collection of all Poisson square-mean almost automorphic stochastic processes $F: \mathbb{R}\times \mathcal{L}^{2}(P,H)\times V\rightarrow \mathcal{L}^{2}(P,H)$ is denoted by $PSAA(\mathbb{R}\times \mathcal{L}^{2}(P,H)\times V, \mathcal{L}^{2}(P,H))$.

In the following lemma, we give some properties of Poisson
  square-mean almost automorphic stochastic processes. \\
 \textbf{Lemma 2.27.}  If $F, F_{1}, F_{2}:\mathbb{R}\times \mathcal{L}^{2}(P,H)\times V\rightarrow \mathcal{L}^{2}(P,H)$ are all Poisson square-mean almost automorphic stochastic processes in $t$ for each $Y\in \mathcal{L}^{2}(P,H)$, then

 (1) $F_{1}+F_{2}$ is Poisson square-mean almost automorphic.

 (2) $\lambda F$ is Poisson square-mean almost automorphic for every scalar $\lambda$.

 (3) For every $Y\in \mathcal{L}^{2}(P,H)$, there exists a constant $M>0$ such that
 \begin{equation*}
 \sup_{t\in \mathbb{R}}\int_{V}\mathbb{E}\|F(t,Y,x)\|^{2}\nu(dx)\leq M.
  \end{equation*}
\textbf{ Proof.} Since (1) and (2) are obvious, we only need to prove (3). If
 \begin{equation*}
 \sup_{t\in \mathbb{R}}\int_{V}\mathbb{E}\|F(t,Y,x)\|^{2}\nu(dx)=\infty
  \end{equation*}
for some $Y\in \mathcal{L}^{2}(P,H)$, then
 there exists a subsequence of real numbers $\{s'_{n}\}$ and a sequence $\{Y'_{n}\}\subset \mathcal{L}^{2}(P,H)$ such that
  \begin{equation*}
 \lim_{n\rightarrow\infty}\int_{V}\mathbb{E}\|F(s'_{n},Y'_{n},x)\|^{2}\nu(dx)=\infty.
 \end{equation*}
 Since $F\in PSAA(\mathbb{R}\times \mathcal{L}^{2}(P,H)\times V, \mathcal{L}^{2}(P,H))$, there exists a subsequence $\{s_{n}\}\subset \{s'_{n}\}$ and a function $\widetilde{F}: \mathbb{R}\times \mathcal{L}^{2}(P,H)\times V\rightarrow\mathcal{L}^{2}(P,H)$ such that
 \begin{equation}\label{pap}
 \lim_{n\rightarrow \infty}\int_{V}\mathbb{E}\|F(t+s_{n},Y,x)-\widetilde{F}(t,Y,x)\|^{2}\nu(dx)=0.
 \end{equation}
 Since $\int_{V}\mathbb{E}\|\widetilde{F}(t,Y,x)\|^{2}\nu(dx)<\infty$ for any fixed $t \in\mathbb{R}$, in particular, taking $t=0$ in (\ref{pap}), we have
 \begin{equation*}
 \lim_{n\rightarrow \infty}\int_{V}\mathbb{E}\|F(s_{n},Y,x)\|^{2}\nu(dx)<\infty,
 \end{equation*}
 a contradiction. \ \ $\Box$\\
\textbf{Lemma 2.28.} Let $F: \mathbb{R}\times \mathcal{L}^{2}(P,H)\times U\rightarrow \mathcal{L}^{2}(P,H)$, $(t,Y,x)\mapsto F(t,Y,x)$ be Poisson square-mean almost automorphic in $t\in \mathbb{R}$ for each $Y\in \mathcal{L}^{2}(P,H)$, and assume that $F$ satisfies the Lipschitz condition in the following sense:
 \begin{equation*}
 \int_{V}\mathbb{E}\|F(t,Y,x)-F(t,Z,x)\|^{2}\nu(dx)\leq L\mathbb{E}\|Y-Z\|^{2}
 \end{equation*}
 for all $Y,Z\in \mathcal{L}^{2}(P,H)$ and each $t\in \mathbb{R}$, where $L>0$ is independent of $t$. Then for any square-mean almost automorphic process $Y: \mathbb{R}\rightarrow \mathcal{L}^{2}(P,H)$, the stochastic process $J: \mathbb{R}\times V\rightarrow \mathcal{L}^{2}(P,H)$ given by $J(t,x):=F(t,Y(t),x)$ is Poisson square-mean almost automorphic. \\
  \textbf{Proof.} Let $\{s'_{n}\}$ be a sequence of real numbers. Since $Y$ is square-mean almost automorphic and $J$ is Poisson square-mean almost automorphic, there exist a subsequence $\{s_{n}\}$ of $\{s'_{n}\}$ such that for some stochastic process $Z: \mathbb{R}\rightarrow \mathcal{L}^{2}(P,H)$ and for each $t\in \mathbb{R}$,
  \begin{equation}\label{right}
  \lim_{n\rightarrow \infty}\mathbb{E}\|Y(t+s_{n})-Z(t)\|^{2}=0,
  \end{equation}
 and for some function $\widetilde{F}: \mathbb{R}\times \mathcal{L}^{2}(P,H)\times V\rightarrow \mathcal{L}^{2}(P,H)$,
  \begin{equation}\label{lim}
  \lim_{n\rightarrow \infty}\int_{V}\mathbb{E}\|F(t+s_{n},Y,x)-\widetilde{F}(t,Y,x)\|^{2}\nu(dx)=0
  \end{equation}
 with any $t\in \mathbb{R}$, $Y\in \mathcal{L}^{2}(P,H)$.

 We consider the function $\widetilde{J}: \mathbb{R}\times V\rightarrow \mathcal{L}^{2}(P,H)$ defined by $\widetilde{J}(t,x):=\widetilde{F}(t,Z(t),x)$. Observe that
\begin{align*}
J(t+s_{n})-\widetilde{J}(t,x)&=F(t+s_{n},Y(t+s_{n}),x)-F(t+s_{n},Z(t),x)\\
&\quad+F(t+s_{n},Z(t),x)-\widetilde{F}(t,Z(t),x),
\end{align*}
then we have
\begin{align}\label{diff}
&\int_{V}\mathbb{E}\|J(t+s_{n},x)-\widetilde{J}(t,x)\|^{2}\nu(dx)\nonumber\\
&\leq2\int_{V}\mathbb{E}\|F(t+s_{n},Y(t+s_{n}),x)-F(t+s_{n},Z(t),x)\|^{2}\nu(dx)\nonumber\\
&\quad+2\int_{V}\mathbb{E}\|F(t+s_{n},Z(t),x)-\widetilde{F}(t,Z(t),x)\|^{2}\nu(dx)\nonumber\\
&\leq2L\mathbb{E}\|Y(t+s_{n})-Z(t)\|^{2}+2\int_{V}\mathbb{E}\|F(t+s_{n},Z(t),x)-\widetilde{F}(t,Z(t),x)\|^{2}\nu(dx).
\end{align}
Put (\ref{right}) and (\ref{lim}) into (\ref{diff}) to get
\begin{align*}
\lim_{n\rightarrow \infty}\int_{V}\mathbb{E}\|J(t+s_{n},x)-\widetilde{J}(t,x)\|^{2}\nu(dx)=0.
\end{align*}
Similarly, we can prove that
\begin{align*}
\lim_{n\rightarrow \infty}\int_{V}\mathbb{E}\|\widetilde{J}(t-s_{n},x)-J(t,x)\|^{2}\nu(dx)=0.
\end{align*}
Hence $J(t,x)$ is Poisson square-mean almost automorphic. \ \ $\Box$ \\
\textbf{Definition 2.29.} A stochastic process $J: \mathbb{R}\times V\rightarrow \mathcal{L}^{2}(P,H)$ is said to be Poisson  square-mean pseudo almost automorphic in $t\in \mathbb{R}$ for each $Y\in \mathcal{L}^{2}(P,H)$ if $J$ is Poisson stochastically continuous and it it can be decomposed as $J=g+\varphi$, where $g\in PSAA(\mathbb{R}\times V,\mathcal{L}^{2}(P,H))$ and $\varphi\in PSBC_{0}(\mathbb{R}\times V, \mathcal{L}^{2}(P,H))$. The collection of all such stochastic processes is denoted by $PSPAA(\mathbb{R}\times V, \mathcal{L}^{2}(P,H))$. \\
\textbf{Definition 2.30.} A stochastic process $J: \mathbb{R}\times V\rightarrow \mathcal{L}^{2}(P,H)$ is said to be Poisson  square-mean weighted pseudo almost automorphic about $\rho\in \mathcal{U}_{\infty}$ in $t\in \mathbb{R}$ if $F$ is Poisson stochastically continuous and it it can be decomposed as $F=g+\varphi$, where $g\in PSAA(\mathbb{R}\times V,\mathcal{L}^{2}(P,H))$ and $\varphi\in PSBC_{0}(\mathbb{R}\times V, \rho)$. The collection of all such stochastic processes is denoted by $PSWPAA(\mathbb{R}\times V, \rho)$. \\

Set
\begin{equation*}
PSAA(\mathbb{R}\times \mathcal{L}^{2}(P,H)\times V, \rho)=\{F(t,Y,x)\in PSAA(\mathbb{R}\times V, \rho)| \ \mbox{for any} \ Y\in \mathcal{L}^{2}(P,H)\},
\end{equation*}
and
\begin{equation*}
PSBC_{0}(\mathbb{R}\times \mathcal{L}^{2}(P,H)\times V, \rho)=\{F(t,Y,x)\in PSBC_{0}(\mathbb{R}\times V, \rho)| \ \mbox{for any} \ Y\in \mathcal{L}^{2}(P,H)\}.
\end{equation*}
\textbf{Definition 2.31.} A stochastic process $F: \mathbb{R}\times \mathcal{L}^{2}(P,H)\times V\rightarrow \mathcal{L}^{2}(P,H)$, $(t,Y,x)\mapsto F(t,Y,x)$ is said to be Poisson  square-mean weighted pseudo almost automorphic about $\rho\in \mathcal{U}_{\infty}$ in $t\in \mathbb{R}$ for each $Y\in \mathcal{L}^{2}(P,H)$ if $F$ Poisson stochastically continuous and it it can be decomposed as $F=g+\varphi$, where $g\in PSAA(\mathbb{R}\times \mathcal{L}^{2}(P,H)\times V,\mathcal{L}^{2}(P,H))$ and $\varphi\in PSBC_{0}(\mathbb{R}\times \mathcal{L}^{2}(P,H)\times V, \rho)$. The collection of all such stochastic processes is denoted by $PSWPAA(\mathbb{R}\times \mathcal{L}^{2}(P,H)\times V, \rho)$. \\
\textbf{Lemma 2.32.}  Assume $J\in PSBC(\mathbb{R}\times V, \mathcal{L}^{2}(P,H))$. Then $J\in PSBC_{0}(\mathbb{R}\times V, \rho)$, where $\rho\in \mathcal{U}_{\infty}$ if and only if for any $\varepsilon>0$,
\begin{equation*}
\lim_{r\rightarrow +\infty}\frac{1}{m(r,\rho)}\int_{M_{r,\varepsilon}(\varphi)}\rho(t)dt=0,
\end{equation*}
where
\begin{equation*}
M_{r,\varepsilon}(\varphi)=\{t\in[-r,r]\big|\int_{V}\mathbb{E}\|J(t,x)\|^{2}\nu(dx)\geq \varepsilon\}.
\end{equation*}

The proof of Lemma 2.32 is included in the Appendix. \\
   \textbf{Theorem 2.33.}  If $\rho\in \mathcal{U}^{\infty}$, $F=g+\varphi\in PSWPAA(\mathbb{R}\times \mathcal{L}^{2}(P,H)\times V, \rho)$ with $g\in PSAA(\mathbb{R}\times \mathcal{L}^{2}(P,H)\times V,\mathcal{L}^{2}(P,H))$ and $\varphi\in PSBC_{0}(\mathbb{R}\times \mathcal{L}^{2}(P,H)\times V, \rho)$. Assume that $F$ and $g$ are Lipschitzian in $Y$ uniformly in $t\in \mathbb{R}$, that is for all $Y,Z\in \mathcal{L}^{2}(P,H)$ and $t\in \mathbb{R}$,
    \begin{align*}
    &\int_{V}\mathbb{E}\|F(t,Y,x)-F(t,Z,x)\|^{2}\nu(dx)\leq LE\|Y-Z\|^{2},\\
     &\int_{V}\mathbb{E}\|g(t,Y,x)-g(t,Z,x)\|^{2}\nu(dx)\leq LE\|Y-Z\|^{2},
    \end{align*}
for some constant $L>0$ independent of $t$. Then for any $Y\in SWPAA (\mathbb{R}, \rho)$, the stochastic process $J: \mathbb{R}\times V\rightarrow \mathcal{L}^{2}(P,H)$ given by $J(t,x):=F(t,Y(t),x)$ is Poisson square-mean weighted  pseudo almost automorphic.

The proof of Theorem 2.33 is included in the Appendix. \\
  \textbf{Definition 2.34.} (see \cite{CZ}.) An  $\mathcal{L}^{2}$-continuous stochastic process $f(t,s): \mathbb{R}\times \mathbb{R} \rightarrow \mathcal{L}^{2}(P,H)$ is said to be square-mean bi-almost automorphic if for every sequence of real numbers $\{s'_{n}\}$, there exists a subsequence $\{s_{n}\}$ and a continuous function $g: \mathbb{R}\times \mathbb{R} \rightarrow \mathcal{L}^{2}(P,H)$ such that
  \begin{equation*}
  \lim_{n\rightarrow \infty}\mathbb{E}\|f(t+s_{n},s+s_{n})-g(t,s)\|^{2}=0
  \end{equation*}
  and
  \begin{equation*}
  \lim_{n\rightarrow \infty}\mathbb{E}\|g(t-s_{n},s-s_{n})-f(t,s)\|^{2}=0
  \end{equation*}
  for each $t,s\in \mathbb{R}$.

  The collection of all square-mean bi-almost automorphic processes is denoted by $SBAA(\mathbb{R}\times \mathbb{R}, \mathcal{L}^{2}(P,H))$.
\section{Square-mean almost automorphic solutions of nonautonomous linear stochastic differential equations}
Consider the following linear nonautonomous stochastic differential equation
\begin{align}\label{levy}
dY(t)=&A(t)Y(t-)dt+f(t)dt+g(t)dW(t)+\int_{|x|_{V}<1}F(t,x)\widetilde{N}(dt,dx)\nonumber\\
&\quad+\int_{|x|_{V}\geq 1}G(t,x)N(dt,dx),
\end{align}
where $f: \mathbb{R}\rightarrow \mathcal{L}^{2}(P,H)$, $g: \mathbb{R}\rightarrow L(V, \mathcal{L}^{2}(P,H))$, $F: \mathbb{R}\times V\rightarrow \mathcal{L}^{2}(P,H))$, $G: \mathbb{R}\times V\rightarrow \mathcal{L}^{2}(P,H))$ are stochastic processes, $W$ and $N$ are the L$\acute{e}$vy-It$\hat{o}$ decomposition components of the two-sided L$\acute{e}$vy process $L$ with assumptions stated in Section 2.1, $A(t)$ satisfies the ``Acquistapace-Terreni" condition \cite{P}, that is:

$(H_{1})$ There exist constants $\lambda_{0}\geq 0$, $\theta\in (\frac{\pi}{2},\pi)$, $K_{1}, K_{2}\geq 0$, and $\alpha_{1},\beta_{1}\in(0,1]$ with
$\alpha_{1}+\beta_{1}>1$  such that
\begin{align*}
\Sigma_{\theta}\cup \{0\}\subset \rho(A(t)-\lambda_{0}), \  \|R(\lambda,A(t)-\lambda_{0})\|\leq \frac{K_{1}}{1+|\lambda|},
\end{align*}
and
\begin{align*}
\|[A(t)-\lambda_{0}]R(\lambda,A(t)-\lambda_{0})[R(\lambda_{0},A(t))-R(\lambda_{0},A(s))]\|\leq K_{2}|t-s|^{\alpha_{1}}|\lambda|^{-\beta_{1}}
\end{align*}
for $t,s\in\mathbb{R}$, $\lambda\in \Sigma_{\theta}:=\{\lambda\in \mathbb{C} \backslash \{0\}: \ |\arg \lambda|\leq \theta\}$.

If the condition $(H_{1})$ is satisfied, then from \cite{P}, it follows that there exists a unique evolution family $\{U(t,s)\}_{-\infty<s\leq t<+\infty}$ on a Banach space $X$, $U$ is the solution of the following equation
\begin{align*}
\frac{d\tau}{dt}=A(t)\tau(t), \ t\geq s, \ \tau(s)=a\in X.
\end{align*}

We further assume that:

$(H_{2})$ $U(t,s)$ is an exponential stable evolution family on $\mathcal{L}^{2}(P,H)$, that is, there exist constants $M,\delta>0$ such that
$\|U(t,s)\|\leq Me^{-\delta(t-s)}$, $t\geq s$.

$(H_{3})$ $U(t,s)\tau\in SBAA(\mathbb{R}\times \mathbb{R}, \mathcal{L}^{2}(P,H))$ uniformly for all $\tau$ in any bounded subset of $\mathcal{L}^{2}(P,H)$.

$(H_{4})$ $U(t,s)y\in SBAA(\mathbb{R}\times \mathbb{R}, \mathcal{L}^{2}(P,H))$ uniformly for all $y\in \mathbb{B}$, where
\begin{align*}
\mathbb{B}:=\{y: \mathbb{R}\times \mathcal{L}^{2}(P,H) \times V\rightarrow \mathcal{L}^{2}(P,H), \ \sup_{t\in \mathbb{R}}\int_{V}\mathbb{E}\|y(t,Y,x)\|^{2}\nu(dx)<\infty\}.
\end{align*}
\textbf{Definition 3.1.} An $\mathcal{F}_{t}$-progressively measurable process $\{Y(t)\}_{t\in \mathbb{R}}$ is called a mild solution of (\ref{levy})
if it satisfies the stochastic integral equation
\begin{align}\label{mild}
Y(t)&=U(t,a)Y(a)+\int_{a}^{t}U(t,s)f(s)ds+\int_{a}^{t}U(t,s)g(s)dW(s)\nonumber\\
&\quad+\int_{a}^{t}\int_{|x|_{V}<1}U(t,s)F(s,x)\widetilde{N}(ds,dx)\nonumber\\
&\quad+\int_{a}^{t}\int_{|x|_{V}\geq 1}U(t,s)G(s,x)N(ds,dx),
\end{align}
for all $t\geq a$ and each $a \in\mathbb{R}$.\\
\textbf{Theorem 3.2.} If assumptions $(H_{1})-(H_{4})$ hold, and $f\in SAA(\mathbb{R}, \mathcal{L}^{2}(P,H))$, $g\in SAA(\mathbb{R}, L(V,\mathcal{L}^{2}(P,H)))$, $F,G\in PSAA(\mathbb{R}\times V, \mathcal{L}^{2}(P,H))$, then the equation (\ref{levy}) has a unique square-mean almost automorphic mild solution. \\
\textbf{Proof.} The process
\begin{align}\label{process}
Y(t)&=\int_{-\infty}^{t}U(t,s)f(s)ds+\int_{-\infty}^{t}U(t,s)g(s)dW(s)\nonumber\\
&\quad+\int_{-\infty}^{t}\int_{|x|_{V}<1}U(t,s)F(s,x)\widetilde{N}(ds,dx)\nonumber\\
&\quad+\int_{-\infty}^{t}\int_{|x|_{V}\geq1}U(t,s)G(s,x)N(ds,dx)
\end{align}
satisfies the equation (\ref{mild}), so it is a mild solution of (\ref{levy}). In the following, we will prove the process (\ref{process}) is the unique square-mean almost automorphic solution of (\ref{levy}) and we divide the proof into three steps.

\emph{Step 1. } The $\mathcal{L}^{2}$-continuity of $Y(t)$. Set $H(t)=\int_{-\infty}^{t}U(t,s)f(s)ds+\int_{-\infty}^{t}U(t,s)g(s)dW(s)$. Similar to
the proof of \cite[Theorem 4.1]{CZ} with minor modifications, it follows that $H(t)$ is $\mathcal{L}^{2}$-continuous. We only need to prove that
\begin{align*}
\int_{-\infty}^{t}\int_{|x|_{V}<1}U(t,s)F(s,x)\widetilde{N}(ds,dx)+\int_{-\infty}^{t}\int_{|x|_{V}\geq1}U(t,s)G(s,x)N(ds,dx)
\end{align*}
is $\mathcal{L}^{2}$-continuous.

Fix $t_{0}\in \mathbb{R}$. Let
$N_{1}(\sigma,x):=N(\sigma+t-t_{0},x)-N(t-t_{0},x)$ and $\widetilde{N}_{1}(\sigma,x):=\widetilde{N}(\sigma+t-t_{0},x)-\widetilde{N}(t-t_{0},x)$ for each $\sigma\in \mathbb{R}$. It is easy to show that $N_{1}$ is also a Poisson random measure and has the same law as $N$, so are $\widetilde{N}_{1}$ and $\widetilde{N}$, moreover, $\widetilde{N}_{1}$  is the compensated Poisson measure of $N_{1}$. Let $\sigma=s-(t-t_{0})$, we have
\begin{align}\label{eleven}
&\mathbb{E}\big\|\int_{-\infty}^{t}\int_{|x|_{V}<1}U(t,s)F(s,x)\widetilde{N}(ds,dx)-\int_{-\infty}^{t_{0}}\int_{|x|_{V}<1}U(t_{0},s)F(s,x)\widetilde{N}(ds,dx)\big\|^{2}\nonumber\\
&\quad+\mathbb{E}\big\|\int_{-\infty}^{t}\int_{|x|_{V}\geq1}U(t,s)G(s,x)N(ds,dx)-\int_{-\infty}^{t_{0}}\int_{|x|_{V}\geq1}U(t_{0},s)G(s,x)N(ds,dx)\big\|^{2}\nonumber\\
&=\mathbb{E}\big\|\int_{-\infty}^{t_{0}}\int_{|x|_{V}<1}U(t,\sigma+t-t_{0})F(\sigma+t-t_{0},x)\widetilde{N}_{1}(d\sigma,dx)\nonumber\\
&\quad-\int_{-\infty}^{t_{0}}\int_{|x|_{V}<1}U(t_{0},s)F(s,x)\widetilde{N}(ds,dx)\big\|^{2}\nonumber\\
&\quad+\mathbb{E}\big\|\int_{-\infty}^{t_{0}}\int_{|x|_{V}\geq1}U(t,\sigma+t-t_{0})G(\sigma+t-t_{0},x)N_{1}(d\sigma,dx)\nonumber\\
&\quad-\int_{-\infty}^{t_{0}}\int_{|x|_{V}\geq1}U(t_{0},s)G(s,x)N(ds,dx)\big\|^{2}\nonumber\\
&=\mathbb{E}\big\|\int_{-\infty}^{t_{0}}\int_{|x|_{V}<1}[U(t,\sigma+t-t_{0})F(\sigma+t-t_{0},x)-U(t_{0},\sigma)F(\sigma,x)]\widetilde{N}_{1}(d\sigma,dx)\big\|^{2}\nonumber\\
&\quad+\mathbb{E}\big\|\int_{-\infty}^{t_{0}}\int_{|x|_{V}\geq1}[U(t,\sigma+t-t_{0})G(\sigma+t-t_{0},x)-U(t_{0},\sigma)G(\sigma,x)]N_{1}(d\sigma,dx)\big\|^{2}\nonumber\\
&\leq \mathbb{E}\big\|\int_{-\infty}^{t_{0}}\int_{|x|_{V}<1}[U(t,\sigma+t-t_{0})F(\sigma+t-t_{0},x)-U(t_{0},\sigma)F(\sigma,x)]\widetilde{N}_{1}(d\sigma,dx)\big\|^{2}\nonumber\\
&\quad+2\mathbb{E}\big\|\int_{-\infty}^{t_{0}}\int_{|x|_{V}\geq1}[U(t,\sigma+t-t_{0})G(\sigma+t-t_{0},x)-U(t_{0},\sigma)G(\sigma,x)]\widetilde{N}_{1}(d\sigma,dx)\big\|^{2}\nonumber\\
&\quad+2\mathbb{E}\big\|\int_{-\infty}^{t_{0}}\int_{|x|_{V}\geq1}[U(t,\sigma+t-t_{0})G(\sigma+t-t_{0},x)-U(t_{0},\sigma)G(\sigma,x)]\nu(dx)d\sigma\big\|^{2}\nonumber\\
&=I_{1}+I_{2}+I_{3}.
\end{align}
For $I_{1}$, by properties of the integral for the Poisson random measure, we have
\begin{align}\label{add}
I_{1}&=\mathbb{E}\big\|\int_{-\infty}^{t_{0}}\int_{|x|_{V}<1}[U(t,\sigma+t-t_{0})F(\sigma+t-t_{0},x)-U(t_{0},\sigma)F(\sigma,x)]\widetilde{N}_{1}(d\sigma,dx)\big\|^{2}\nonumber\\
&\leq \int_{-\infty}^{t_{0}}\int_{|x|_{V}<1}\mathbb{E}\|U(t,\sigma+t-t_{0})F(\sigma+t-t_{0},x)-U(t_{0},\sigma)F(\sigma,x)\|^{2}\nu(dx)d\sigma\nonumber\\
&=\int_{-\infty}^{t_{0}}\int_{|x|_{V}<1}\mathbb{E}\|U(t,\sigma+t-t_{0})(F(\sigma+t-t_{0},x)-F(\sigma,x))\nonumber\\
&\quad+(U(t,\sigma+t-t_{0})F(\sigma,x)-U(t_{0},\sigma)F(\sigma,x))\|^{2}\nu(dx)d\sigma\nonumber\\
&\leq 2\int_{-\infty}^{t_{0}}\int_{|x|_{V}<1}\|U(t,\sigma+t-t_{0})\|^{2}\mathbb{E}\|F(\sigma+t-t_{0},x)-F(\sigma,x)\|^{2}\nu(dx)d\sigma\nonumber\\
&\quad+2\int_{-\infty}^{t_{0}}\int_{|x|_{V}<1}\mathbb{E}\|U(t,\sigma+t-t_{0})F(\sigma,x)-U(t_{0},\sigma)F(\sigma,x)\|^{2}\nu(dx)d\sigma.
\end{align}
By the exponential dissipation property $(H_{2})$ of $U(t,s)$, we get
\begin{align}\label{u}
&\int_{|x|_{V}<1}\|U(t,\sigma+t-t_{0})\|^{2}\mathbb{E}\|F(\sigma+t-t_{0},x)-F(\sigma,x)\|^{2}\nu(dx)\nonumber\\
&\leq 2\|U(t,\sigma+t-t_{0})\|^{2}\left(\int_{|x|_{V}<1}\mathbb{E}\|F(\sigma+t-t_{0},x)\|^{2}\nu(dx)+\int_{|x|_{U}<1}E\|F(\sigma,x)\|^{2}\nu(dx)\right)\nonumber\\
&\leq 4M^{2}e^{-2\delta(t_{0}-\sigma)}\sup_{\sigma\in \mathbb{R}}\int_{|x|_{V}<1}\mathbb{E}\|F(\sigma,x)\|^{2}\nu(dx).
\end{align}
By (3) of Lemma 2.27,
\begin{align}\label{t}
\int_{-\infty}^{t_{0}}4M^{2}e^{-2\delta(t_{0}-\sigma)}\sup_{\sigma\in \mathbb{R}}\int_{|x|_{V}<1}\mathbb{E}\|F(\sigma,x)\|^{2}\nu(dx)d\sigma<\infty.
\end{align}
By the continuity of Poisson square-mean almost automorphic functions,
\begin{align}\label{c}
\int_{|x|_{V}<1}\mathbb{E}\|F(\sigma+t-t_{0},x)-F(\sigma,x)\|^{2}\nu(dx)\rightarrow 0 \ \ \mbox{as} \ t\rightarrow t_{0}.
\end{align}
From (\ref{u}), (\ref{t}), (\ref{c}) and the Lebesgue dominated convergence theorem, it follows that
\begin{align}\label{Lebesgue}
\int_{-\infty}^{t_{0}}\int_{|x|_{V}<1}\|U(t,\sigma+t-t_{0})\|^{2}\mathbb{E}\|F(\sigma+t-t_{0},x)-F(\sigma,x)\|^{2}\nu(dx)d\sigma\rightarrow 0 \ \ \mbox{as} \ t\rightarrow t_{0}.
\end{align}
Thanks to the strong continuity of $U(t,s)$,
\begin{align*}
E\big\|U(t,\sigma+t-t_{0})F(\sigma,x)-U(t_{0},\sigma)F(\sigma,x)\big\|^{2}\rightarrow 0 \ \ \mbox{as} \ t\rightarrow t_{0}.
\end{align*}
By the exponential dissipation property of $U(t,s)$, we have
\begin{align*}
&\int_{|x|_{V}<1}\mathbb{E}\|U(t,\sigma+t-t_{0})F(\sigma,x)-U(t_{0},\sigma)F(\sigma,x)\|^{2}\nu(dx)\nonumber\\
&\leq 4M^{2}e^{-2\delta(t_{0}-\sigma)}\sup_{\sigma\in \mathbb{R}}\int_{|x|_{V}<1}\mathbb{E}\|F(\sigma,x)\|^{2}\nu(dx).
\end{align*}
By the Lebesgue dominated convergence theorem,
\begin{align}\label{dissipation}
\int_{-\infty}^{t_{0}}\int_{|x|_{V}<1}\mathbb{E}\|U(t,\sigma+t-t_{0})F(\sigma,x)-U(t_{0},\sigma)F(\sigma,x)\|^{2}\nu(dx)d\sigma\rightarrow 0 \ \ \mbox{as} \ t\rightarrow t_{0}.
\end{align}
By (\ref{add}), (\ref{Lebesgue}) and (\ref{dissipation}), we have
\begin{align}\label{I1}
\lim_{t\rightarrow t_{0}}I_{1}=0.
\end{align}
Similarly,
\begin{align}\label{I2}
\lim_{t\rightarrow t_{0}}I_{2}=0.
\end{align}
For $I_{3}$, we have
\begin{align}\label{dx}
&\mathbb{E}\big\|\int_{-\infty}^{t_{0}}\int_{|x|_{V}\geq 1}[U(t,\sigma+t-t_{0})G(\sigma+t-t_{0},x)-U(t_{0},\sigma)G(\sigma,x)]\nu(dx)d\sigma\big\|^{2}\nonumber\\
&=\mathbb{E}\big\|\int_{-\infty}^{t_{0}}\int_{|x|_{V}\geq 1}[U(t,\sigma+t-t_{0})(G(\sigma+t-t_{0},x)-G(\sigma,x))\nonumber\\
&\quad+(U(t,\sigma+t-t_{0})G(\sigma,x)-U(t_{0},\sigma)G(\sigma,x))]\nu(dx)d\sigma\big\|^{2}\nonumber\\
&\leq 2\mathbb{E}\big\|\int_{-\infty}^{t_{0}}\int_{|x|_{V}\geq 1}U(t,\sigma+t-t_{0})[G(\sigma+t-t_{0},x)-G(\sigma,x)]\nu(dx)d\sigma\big\|^{2}\nonumber\\
&\quad+2\mathbb{E}\big\|\int_{-\infty}^{t_{0}}\int_{|x|_{V}\geq 1}[U(t,\sigma+t-t_{0})G(\sigma,x)-U(t_{0},\sigma)G(\sigma,x)]\nu(dx)d\sigma\big\|^{2}.
\end{align}
By the Cauchy-Schwarz inequality, we have
\begin{align}\label{Schwarz}
&\mathbb{E}\big\|\int_{-\infty}^{t_{0}}\int_{|x|_{V}\geq 1}U(t,\sigma+t-t_{0})[G(\sigma+t-t_{0},x)-G(\sigma,x)]\nu(dx)d\sigma\big\|^{2}\nonumber\\
&\leq M^{2}\mathbb{E}\left(\int_{-\infty}^{t_{0}}\int_{|x|_{V}\geq 1}e^{-\frac{\delta(t_{0}-\sigma)}{2}}\cdot e^{-\frac{\delta(t_{0}-\sigma)}{2}}\cdot \|G(\sigma+t-t_{0},x)-G(\sigma,x)\|\nu(dx)d\sigma \right)^{2}\nonumber\\
&\leq M^{2}\int_{-\infty}^{t_{0}}e^{-\delta(t_{0}-\sigma)}d\sigma\int_{|x|_{V}\geq 1}\nu(dx)\int_{-\infty}^{t_{0}}\int_{|x|_{V}\geq 1}e^{-\delta(t_{0}-\sigma)}\mathbb{E}\|G(\sigma+t-t_{0},x)-G(\sigma,x)\|^{2}\nu(dx)d\sigma\nonumber\\
&=\frac{M^{2}c}{\delta}\int_{-\infty}^{t_{0}}\int_{|x|_{U}\geq 1}e^{-\delta(t_{0}-\sigma)}\mathbb{E}\|G(\sigma+t-t_{0},x)-G(\sigma,x)\|^{2}\nu(dx)d\sigma.
\end{align}
From (\ref{measure}), it follows that $c<\infty$. By the Lebesgue dominated convergence theorem and the similar argument for (\ref{Lebesgue}), we obtain
\begin{align}\label{dominated}
\mathbb{E}\big\|\int_{-\infty}^{t_{0}}\int_{|x|_{V}\geq 1}U(t,\sigma+t-t_{0})[G(\sigma+t-t_{0},x)-G(\sigma,x)]\nu(dx)d\sigma\big\|^{2}\rightarrow 0\ \ \mbox{as} \ t\rightarrow t_{0}.
\end{align}
From the strong continuity of $U(t,s)$, it follows that
\begin{align*}
\mathbb{E}\big\|U(t,\sigma+t-t_{0})G(\sigma,x)-U(t_{0},\sigma)G(\sigma,x)\big\|^{2}\rightarrow 0 \ \ \mbox{as} \ t\rightarrow t_{0}.
\end{align*}
By the exponential dissipation property of $U(t,s)$, we have
\begin{align*}
&\int_{|x|_{V}\geq1}\mathbb{E}\|U(t,\sigma+t-t_{0})G(\sigma,x)-U(t_{0},\sigma)G(\sigma,x)\|^{2}\nu(dx)\nonumber\\
&\leq 4M^{2}e^{-2\delta(t_{0}-\sigma)}\sup_{\sigma\in \mathbb{R}}\int_{|x|_{V}\geq1}\mathbb{E}\|G(\sigma,x)\|^{2}\nu(dx).
\end{align*}
By (3) of Lemma 2.27,
\begin{align*}
\int_{-\infty}^{t_{0}}4M^{2}e^{-2\delta(t_{0}-\sigma)}\sup_{\sigma\in \mathbb{R}}\int_{|x|_{V}\geq1}\mathbb{E}\|G(\sigma,x)\|^{2}\nu(dx)d\sigma<\infty.
\end{align*}
By the Lebesgue dominated convergence theorem, we have
\begin{align*}
\int_{-\infty}^{t_{0}}\int_{|x|_{V}\geq 1}[U(t,\sigma+t-t_{0})G(\sigma,x)-U(t_{0},\sigma)G(\sigma,x)]\nu(dx)d\sigma\rightarrow 0
\end{align*}
in $\mathcal{L}^{2}(P,H)$ as $t\rightarrow t_{0}$, that is,
\begin{align}\label{U}
\mathbb{E}\big\|\int_{-\infty}^{t_{0}}\int_{|x|_{V}\geq 1}[U(t,\sigma+t-t_{0})G(\sigma,x)-U(t_{0},\sigma)G(\sigma,x)]\nu(dx)d\sigma\big\|^{2}\rightarrow 0 \ \ \mbox{as} \ t\rightarrow t_{0}.
\end{align}
By (\ref{dx}), (\ref{dominated}) and (\ref{U}), we have
\begin{align}\label{I3}
\lim_{t\rightarrow t_{0}}I_{3}=0.
\end{align}
By (\ref{eleven}), (\ref{I1}),  (\ref{I2}) and (\ref{I3}), we obtain
\begin{align}\label{I}
&\mathbb{E}\big\|\int_{-\infty}^{t}\int_{|x|_{V}<1}U(t,s)F(s,x)\widetilde{N}(ds,dx)-\int_{-\infty}^{t_{0}}\int_{|x|_{V}<1}U(t_{0},s)F(s,x)\widetilde{N}(ds,dx)\big\|^{2}\nonumber\\
&\quad+\mathbb{E}\big\|\int_{-\infty}^{t}\int_{|x|_{V}\geq1}U(t,s)G(s,x)N(ds,dx)-\int_{-\infty}^{t_{0}}\int_{|x|_{V}\geq1}U(t_{0},s)G(s,x)N(ds,dx)\big\|^{2}\rightarrow 0 \nonumber\\
&\mbox{as} \ t\rightarrow t_{0}.
\end{align}
Hence,
\begin{align*}
\int_{-\infty}^{t}\int_{|x|_{V}<1}U(t,s)F(s,x)\widetilde{N}(ds,dx)+\int_{-\infty}^{t}\int_{|x|_{V}\geq1}U(t,s)G(s,x)N(ds,dx)
\end{align*}
is $\mathcal{L}^{2}$-continuous. Therefore, $Y(t)$ is $\mathcal{L}^{2}$-continuous.

\emph{Step 2. Existence.} Similar to
the proof of \cite[Theorem 3.1]{CZ} with minor modifications, it follows that
\begin{align*}
\int_{-\infty}^{t}U(t,s)f(s)ds \ \ \mbox{and} \ \int_{-\infty}^{t}U(t,s)g(s)dW(s)
\end{align*}
are square-mean automorphic.

Let $\{s'_{n}\}$ be an arbitrary sequence of real numbers. Since $F,G\in PSAA(\mathbb{R}\times U, \mathcal{L}^{2}(P,H))$, there exists a subsequence
 $\{s_{n}\}$ of  $\{s'_{n}\}$ and continuous functions $\widetilde{F}, \widetilde{G}$ such that
\begin{align}\label{PSAA}
 \lim_{n\rightarrow \infty}\int_{|x|_{V}<1}\mathbb{E}\|F(t+s_{n},x)-\widetilde{F}(t,x)\|^{2}\nu(dx)=0
\end{align}
and
\begin{align}\label{A}
 \lim_{n\rightarrow \infty}\int_{|x|_{V}\geq1}\mathbb{E}\|G(t+s_{n},x)-\widetilde{G}(t,x)\|^{2}\nu(dx)=0
\end{align}
for each $t\in \mathbb{R}$. \\
By $(H_{4})$, there exists a evolution family $U_{1}(t,s)$ such that
\begin{align}\label{bound}
\lim_{n\rightarrow \infty}\mathbb{E}\|U(t+s_{n},s+s_{n})y-U_{1}(t,s)y\|^{2}=0
\end{align}
and
\begin{align}\label{minus}
\lim_{n\rightarrow \infty}\mathbb{E}\|U_{1}(t-s_{n},s-s_{n})y-U(t,s)y\|^{2}=0
\end{align}
for each $y\in \mathbb{B}$.

Let $N_{1}(\sigma,x):=N(\sigma+s_{n},x)-N(s_{n},x)$ and $\widetilde{N}_{1}(\sigma,x):=\widetilde{N}(\sigma+s_{n},x)-\widetilde{N}(s_{n},x)$ for each $\sigma\in \mathbb{R}$. As in step 1, $N_{1}$ has the same law as $N$ with the compensated Poisson random measure $\widetilde{N}_{1}$. Let $\sigma=s-s_{n}$. By properties of the integral for the Poisson random measure, we have
\begin{align}\label{expect}
&\mathbb{E}\big\|\int_{-\infty}^{t+s_{n}}\int_{|x|_{V}<1}U(t+s_{n},s)F(s,x)\widetilde{N}(ds,dx)-\int_{-\infty}^{t}\int_{|x|_{V}<1}U_{1}(t,s)\widetilde{F}(s,x)\widetilde{N}(ds,dx)\big\|^{2}\nonumber\\
&=\mathbb{E}\big\|\int_{-\infty}^{t}\int_{|x|_{V}<1}U(t+s_{n},\sigma+s_{n})F(\sigma+s_{n},x)\widetilde{N}_{1}(d\sigma,dx)\nonumber\\
&\quad-\int_{-\infty}^{t}\int_{|x|_{V}<1}U_{1}(t,\sigma)\widetilde{F}(\sigma,x)\widetilde{N}_{1}(d\sigma,dx)\big\|^{2}\nonumber\\
&=\mathbb{E}\big\|\int_{-\infty}^{t}\int_{|x|_{V}<1}U(t+s_{n},\sigma+s_{n})[F(\sigma+s_{n},x)-\widetilde{F}(\sigma,x)]\widetilde{N}_{1}(d\sigma,dx)\nonumber\\
&\quad-\int_{-\infty}^{t}\int_{|x|_{V}<1}[U(t+s_{n},\sigma+s_{n})-U_{1}(t,\sigma)]\widetilde{F}(\sigma,x)\widetilde{N}_{1}(d\sigma,dx)\big\|^{2}\nonumber\\
&\leq 2\mathbb{E}\big\|\int_{-\infty}^{t}\int_{|x|_{V}<1}U(t+s_{n},\sigma+s_{n})[F(\sigma+s_{n})-\widetilde{F}(\sigma,x)]\widetilde{N}_{1}(d\sigma,dx)\big\|^{2}\nonumber\\
&\quad+2\mathbb{E}\big\|\int_{-\infty}^{t}\int_{|x|_{V}<1}[U(t+s_{n},\sigma+s_{n})-U_{1}(t,\sigma)]\widetilde{F}(\sigma,x)\widetilde{N}_{1}(d\sigma,dx)\big\|^{2}\nonumber\\
&\leq 2\int_{-\infty}^{t}\int_{|x|_{V}<1}\|U(t+s_{n},\sigma+s_{n})\|^{2}\mathbb{E}\|F(\sigma+s_{n},x)-\widetilde{F}(\sigma,x)\|^{2}\nu(dx)d\sigma\nonumber\\
&\quad+2\int_{-\infty}^{t}\int_{|x|_{V}<1}\mathbb{E}\|[U(t+s_{n},\sigma+s_{n})-U_{1}(t,\sigma)]\widetilde{F}(\sigma,x)\|^{2}\nu(dx)d\sigma\nonumber\\
&\leq 2\int_{-\infty}^{t}M^{2}e^{-2\delta(t-\sigma)}\int_{|x|_{V}<1}\mathbb{E}\|F(\sigma+s_{n},x)-\widetilde{F}(\sigma,x)\|^{2}\nu(dx)d\sigma\nonumber\\
&\quad+2\int_{-\infty}^{t}\int_{|x|_{V}<1}\mathbb{E}\|[U(t+s_{n},\sigma+s_{n})-U_{1}(t,\sigma)]\widetilde{F}(\sigma,x)\|^{2}\nu(dx)d\sigma.
\end{align}
By (\ref{PSAA}) and the Lebesgue dominated convergence theorem, we have
\begin{align}\label{almost}
\lim_{n\rightarrow\infty}\int_{-\infty}^{t}M^{2}e^{-2\delta(t-\sigma)}\int_{|x|_{V}<1}\mathbb{E}\|F(\sigma+s_{n},x)-\widetilde{F}(\sigma,x)\|^{2}\nu(dx)d\sigma=0.
\end{align}
 By (\ref{bound}) and the exponential dissipation property of $U(t,s)$, we get
\begin{align}\label{U1}
\mathbb{E}\|U_{1}(t,s)y\|^{2}\leq 2M^{2}e^{-2\delta(t-s)}\mathbb{E}\|y\|^{2}
\end{align}
for all $t\geq s$ and $y\in \mathbb{B}$.
By (\ref{U1}), we have
\begin{align*}
&\int_{|x|_{V}<1}\mathbb{E}\|[U(t+s_{n},\sigma+s_{n})-U_{1}(t,\sigma)]\widetilde{F}(\sigma,x)\|^{2}\nu(dx)\nonumber\\
&\leq 6M^{2}e^{-2\delta(t-\sigma)}\sup_{\sigma\in \mathbb{R}}\int_{|x|_{V}<1}\mathbb{E}\|\widetilde{F}(\sigma,x)\|^{2}\nu(dx).
\end{align*}
From (3) of Lemma 2.27 and (\ref{PSAA}), it follows that $\widetilde{F}\in \mathbb{B}$. This implies
\begin{align*}
\int_{-\infty}^{t}6M^{2}e^{-2\delta(t-\sigma)}\sup_{\sigma\in \mathbb{R}}\int_{|x|_{V}<1}\mathbb{E}\|\widetilde{F}(\sigma,x)\|^{2}\nu(dx)d\sigma<\infty.
\end{align*}
By (\ref{bound}) and the Lebesgue dominated convergence theorem, we have
\begin{align}\label{converge}
\lim_{n\rightarrow \infty}\int_{-\infty}^{t}\int_{|x|_{V}<1}\mathbb{E}\|[U(t+s_{n},\sigma+s_{n})-U_{1}(t,\sigma)]\widetilde{F}(\sigma,x)\|^{2}\nu(dx)d\sigma= 0.
\end{align}
From (\ref{expect}), (\ref{almost}) and (\ref{converge}), it follows that
\begin{align*}
\lim_{n\rightarrow \infty}\mathbb{E}\big\|\int_{-\infty}^{t+s_{n}}\int_{|x|_{V}<1}U(t+s_{n},s)F(s,x)\widetilde{N}(ds,dx)-\int_{-\infty}^{t}\int_{|x|_{V}<1}U_{1}(t,s)\widetilde{F}(s,x)\widetilde{N}(ds,dx)\big\|^{2}=0.
\end{align*}
We can use the similar step to prove that
\begin{align*}
\lim_{n\rightarrow \infty}\mathbb{E}\big\|\int_{-\infty}^{t-s_{n}}\int_{|x|_{V}<1}U_{1}(t-s_{n},s)F(s,x)\widetilde{N}(ds,dx)-\int_{-\infty}^{t}\int_{|x|_{V}<1}U(t,s)\widetilde{F}(s,x)\widetilde{N}(ds,dx)\big\|^{2}=0.
\end{align*}
Therefore,
\begin{align*}
\int_{-\infty}^{t}\int_{|x|_{U}<1}U(t,s)F(s,x)\widetilde{N}(ds,dx)
\end{align*}
is square-mean automorphic.

By properties of the integral for the Poisson random measure, we have
\begin{align}\label{random}
&\mathbb{E}\big\|\int_{-\infty}^{t+s_{n}}\int_{|x|_{V}\geq1}U(t+s_{n},s)G(s,x)N(ds,dx)-\int_{-\infty}^{t}\int_{|x|_{V}\geq1}U_{1}(t,s)\widetilde{G}(s,x)N(ds,dx)\big\|^{2}\nonumber\\
&=\mathbb{E}\big\|\int_{-\infty}^{t}\int_{|x|_{V}\geq1}U(t+s_{n},\sigma+s_{n})G(\sigma+s_{n},x)N_{1}(d\sigma,dx)\nonumber\\
&\quad-\int_{-\infty}^{t}\int_{|x|_{V}\geq1}U_{1}(t,\sigma)\widetilde{G}(\sigma,x)N_{1}(d\sigma,dx)\big\|^{2}\nonumber\\
&=\mathbb{E}\big\|\int_{-\infty}^{t}\int_{|x|_{V}\geq1}[U(t+s_{n},\sigma+s_{n})G(\sigma+s_{n},x)-U_{1}(t,\sigma)\widetilde{G}(\sigma,x)]N_{1}(d\sigma,dx)\big\|^{2}\nonumber\\
&\leq2\mathbb{E}\big\|\int_{-\infty}^{t}\int_{|x|_{V}\geq1}[U(t+s_{n},\sigma+s_{n})G(\sigma+s_{n},x)-U_{1}(t,\sigma)\widetilde{G}(\sigma,x)]\widetilde{N}_{1}(d\sigma,dx)\big\|^{2}\nonumber\\
&\quad+2\mathbb{E}\big\|\int_{-\infty}^{t}\int_{|x|_{V}\geq1}[U(t+s_{n},\sigma+s_{n})G(\sigma+s_{n},x)-U_{1}(t,\sigma)\widetilde{G}(\sigma,x)]\nu(dx)d\sigma\big\|^{2}\nonumber\\
&\leq 2\int_{-\infty}^{t}\int_{|x|_{V}\geq1}\mathbb{E}\|U(t+s_{n},\sigma+s_{n})G(\sigma+s_{n},x)-U_{1}(t,\sigma)\widetilde{G}(\sigma,x)\|^{2}\nu(dx)d\sigma\nonumber\\
&\quad+2\mathbb{E}\left(\int_{-\infty}^{t}\int_{|x|_{V}\geq1}\|U(t+s_{n},\sigma+s_{n})G(\sigma+s_{n},x)-U_{1}(t,\sigma)\widetilde{G}(\sigma,x)\|\nu(dx)d\sigma\right)^{2}\nonumber\\
&\leq 2\int_{-\infty}^{t}\int_{|x|_{V}\geq1}\mathbb{E}\|U(t+s_{n},\sigma+s_{n})[G(\sigma+s_{n},x)-\widetilde{G}(\sigma,x)]\nonumber\\
&\quad+[U(t+s_{n},\sigma+s_{n})-U_{1}(t,\sigma)]\widetilde{G}(\sigma,x)\|^{2}\nu(dx)d\sigma\nonumber\\
&\quad+2\mathbb{E}(\int_{-\infty}^{t}\int_{|x|_{V}\geq1}\|U(t+s_{n},\sigma+s_{n})[G(\sigma+s_{n},x)-\widetilde{G}(\sigma,x)]\nonumber\\
&\quad+[U(t+s_{n},\sigma+s_{n})-U_{1}(t,\sigma)]\widetilde{G}(\sigma,x)\|\nu(dx)d\sigma)^{2}\nonumber\\
&\leq 4\int_{-\infty}^{t}\int_{|x|_{V}\geq1}\mathbb{E}\|U(t+s_{n},\sigma+s_{n})[G(\sigma+s_{n},x)-\widetilde{G}(\sigma,x)]\|^{2}\nu(dx)d\sigma\nonumber\\
&\quad+4\int_{-\infty}^{t}\int_{|x|_{V}\geq1}\mathbb{E}\|[U(t+s_{n},\sigma+s_{n})-U_{1}(t,\sigma)]\widetilde{G}(\sigma,x)\|^{2}\nu(dx)d\sigma\nonumber\\
&\quad+4\mathbb{E}\left(\int_{-\infty}^{t}\int_{|x|_{V}\geq1}\|U(t+s_{n},\sigma+s_{n})[G(\sigma+s_{n})-\widetilde{G}(\sigma,x)]\|\nu(dx)d\sigma\right)^{2}\nonumber\\
&\quad+4\mathbb{E}\left(\int_{-\infty}^{t}\int_{|x|_{V}\geq1}\|[U(t+s_{n},\sigma+s_{n})-U_{1}(t,\sigma)]\widetilde{G}(\sigma,x)\|\nu(dx)d\sigma\right)^{2}.
\end{align}
Similar to the arguments for (\ref{almost}) and (\ref{converge}), by (\ref{A}), (\ref{bound}) and the Lebesgue dominated convergence theorem, we have \begin{align}\label{lim1}
\lim_{n\rightarrow \infty}\int_{-\infty}^{t}\int_{|x|_{V}\geq1}\mathbb{E}\|U(t+s_{n},\sigma+s_{n})[G(\sigma+s_{n},x)-\widetilde{G}(\sigma,x)]\|^{2}\nu(dx)d\sigma=0.
\end{align}
\begin{align}\label{lim2}
\lim_{n\rightarrow \infty}\int_{-\infty}^{t}\int_{|x|_{V}\geq1}\mathbb{E}\|[U(t+s_{n},\sigma+s_{n})-U_{1}(t,\sigma)]\widetilde{G}(\sigma,x)\|^{2}\nu(dx)d\sigma=0.
\end{align}
From the exponential dissipation property of $U(t,s)$ and the Cauchy-Schwarz inequality, it follows that
\begin{align}\label{lim3}
&\mathbb{E}\left(\int_{-\infty}^{t}\int_{|x|_{V}\geq1}\|U(t+s_{n},\sigma+s_{n})[G(\sigma+s_{n},x)-\widetilde{G}(\sigma,x)]\|\nu(dx)d\sigma\right)^{2}\nonumber\\
&\leq M^{2}\mathbb{E}\left(\int_{-\infty}^{t}\int_{|x|_{V}\geq1}e^{-\frac{\delta(t-\sigma)}{2}}\cdot e^{-\frac{\delta(t-\sigma)}{2}}\cdot\|G(\sigma+s_{n},x)-\widetilde{G}(\sigma,x)\|\nu(dx)d\sigma\right)^{2}\nonumber\\
&\leq M^{2}\int_{-\infty}^{t}e^{-\delta(t-\sigma)}d\sigma\int_{|x|_{V}\geq 1}\nu(dx)\int_{-\infty}^{t}\int_{|x|_{V}\geq 1}e^{-\delta(t-\sigma)}\mathbb{E}\|G(\sigma+s_{n},x)-\widetilde{G}(\sigma,x)\|^{2}\nu(dx)d\sigma\nonumber\\
&\leq \frac{M^{2}}{\delta}\int_{|x|_{V}\geq 1}\nu(dx)\int_{-\infty}^{t}\int_{|x|_{V}\geq 1}e^{-\delta(t-\sigma)}\mathbb{E}\|G(\sigma+s_{n},x)-\widetilde{G}(\sigma,x)\|^{2}\nu(dx)d\sigma.
\end{align}
By (\ref{A}), (\ref{lim3}) and the Lebesgue dominated convergence theorem, we have
\begin{align}\label{lim4}
\lim_{n\rightarrow\infty}\mathbb{E}\left(\int_{-\infty}^{t}\int_{|x|_{V}\geq1}\|U(t+s_{n},\sigma+s_{n})[G(\sigma+s_{n},x)-\widetilde{G}(\sigma,x)]\|\nu(dx)d\sigma\right)^{2}=0.
\end{align}
Let $b\in(0,2\delta)$, by the Cauchy-Schwarz inequality, we have
\begin{align}\label{lim5}
&\mathbb{E}\left(\int_{-\infty}^{t}\int_{|x|_{V}\geq1}\|[U(t+s_{n},\sigma+s_{n})-U_{1}(t,\sigma)]\widetilde{G}(\sigma,x)\|\nu(dx)d\sigma\right)^{2}\nonumber\\
&\leq \mathbb{E}\left( \int_{-\infty}^{t}\int_{|x|_{V}\geq1}e^{-\frac{b(t-\sigma)}{2}}\cdot e^{\frac{b(t-\sigma)}{2}}\cdot \|[U(t+s_{n},\sigma+s_{n})-U_{1}(t,\sigma)]\widetilde{G}(\sigma,x)\|\nu(dx)d\sigma\right)^{2}\nonumber\\
&\leq \int_{-\infty}^{t}e^{-b(t-\sigma)}d\sigma\int_{|x|_{V}\geq 1}\nu(dx)\int_{-\infty}^{t}\int_{|x|_{V}\geq1}e^{b(t-\sigma)}\mathbb{E}\|[U(t+s_{n},\sigma+s_{n})-U_{1}(t,\sigma)]\widetilde{G}(\sigma,x)\|^{2}\nu(dx)d\sigma\nonumber\\
&\leq \frac{c}{b}\int_{-\infty}^{t}\int_{|x|_{V}\geq1}e^{b(t-\sigma)}\mathbb{E}\|[U(t+s_{n},\sigma+s_{n})-U_{1}(t,\sigma)]\widetilde{G}(\sigma,x)\|^{2}\nu(dx)d\sigma.
\end{align}
From (3) of Lemma 2.27 and (\ref{A}), it follows that $\widetilde{G}\in \mathbb{B}$.
By the exponential dissipation of $U(t,s)$ and (\ref{U1}), we get
\begin{align*}
&\int_{|x|_{V}\geq 1}e^{b(t-\sigma)}\mathbb{E}\|[U(t+s_{n},\sigma+s_{n})-U_{1}(t,\sigma)]\widetilde{G}(\sigma,x)\|^{2}\nu(dx)\nonumber\\
&\leq \int_{|x|_{V}\geq 1}e^{b(t-\sigma)}(2\mathbb{E}\|U(t+s_{n},\sigma+s_{n})\widetilde{G}(\sigma,x)\|^{2}+2\mathbb{E}\|U_{1}(t,\sigma)\widetilde{G}(\sigma,x)\|^{2})\nu(dx)\nonumber\\
&\leq 6M^{2}e^{(b-2\delta)(t-\sigma)}\sup_{\sigma\in \mathbb{R}}\int_{|x|_{V}\geq 1}\mathbb{E}\|\widetilde{G}(\sigma,x)\|^{2}\nu(dx).
\end{align*}
Since $b<2\delta$,
\begin{align*}
\int_{-\infty}^{t}6M^{2}e^{(b-2\delta)(t-\sigma)}\sup_{\sigma\in \mathbb{R}}\int_{|x|_{V}\geq 1}\mathbb{E}\|\widetilde{G}(\sigma,x)\|^{2}\nu(dx)d\sigma<\infty.
\end{align*}
By (\ref{bound}),  we have
\begin{align}\label{integral}
\lim_{n\rightarrow \infty}\int_{|x|_{V}\geq1}e^{b(t-\sigma)}\mathbb{E}\|[U(t+s_{n},\sigma+s_{n})-U_{1}(t,\sigma)]\widetilde{G}(\sigma,x)\|^{2}\nu(dx)=0.
\end{align}
By (\ref{lim5}), (\ref{integral}) and the Lebesgue dominated convergence theorem, we have
\begin{align}\label{lim6}
\lim_{n\rightarrow \infty}\mathbb{E}\left(\int_{-\infty}^{t}\int_{|x|_{V}\geq1}\|[U(t+s_{n},\sigma+s_{n})-U_{1}(t,\sigma)]\widetilde{G}(\sigma,x)\|\nu(dx)d\sigma\right)^{2}=0.
\end{align}
By (\ref{random}), (\ref{lim1}), (\ref{lim2}), (\ref{lim4}) and (\ref{lim6}), we obtain
\begin{align*}
\lim_{n\rightarrow \infty}\mathbb{E}\big\|\int_{-\infty}^{t+s_{n}}\int_{|x|_{V}\geq1}U(t+s_{n},s)G(s,x)N(ds,dx)-\int_{-\infty}^{t}\int_{|x|_{V}\geq1}U_{1}(t,s)\widetilde{G}(s,x)N(ds,dx)\big\|^{2}=0.
\end{align*}
We can use the similar step to prove that
\begin{align*}
\lim_{n\rightarrow \infty}\mathbb{E}\big\|\int_{-\infty}^{t-s_{n}}\int_{|x|_{V}\geq1}U_{1}(t-s_{n},s)G(s,x)N(ds,dx)-\int_{-\infty}^{t}\int_{|x|_{V}\geq1}U(t,s)\widetilde{G}(s,x)N(ds,dx)\big\|^{2}=0.
\end{align*}
Therefore,
\begin{align*}
\int_{-\infty}^{t}\int_{|x|_{V}\geq 1}U(t,s)G(s,x)N(ds,dx)
\end{align*}
is square-mean almost automorphic.
So, the process $Y(t)$ is square-mean almost automorphic. \\
\emph{Step 3. Uniqueness.} Assume that $Y(t)$ and $Z(t)$ are both square-mean almost automorphic solutions of (\ref{levy}) with differential initial value $Y(a)$ and $Z(a)$ for some $a\in \mathbb{R}$. Then for $t\geq a$,
\begin{align*}
Y(t)&=U(t,a)Y(a)+\int_{a}^{t}U(t,s)f(s)ds+\int_{a}^{t}U(t,s)g(s)dW(s)\nonumber\\
&\quad+\int_{a}^{t}\int_{|x|_{V}<1}U(t,s)F(s,x)\widetilde{N}(ds,dx)\nonumber\\
&\quad+\int_{a}^{t}\int_{|x|_{V}\geq 1}U(t,s)G(s,x)N(ds,dx),
\end{align*}
\begin{align*}
Z(t)&=U(t,a)Z(a)+\int_{a}^{t}U(t,s)f(s)ds+\int_{a}^{t}U(t,s)g(s)dW(s)\nonumber\\
&\quad+\int_{a}^{t}\int_{|x|_{V}<1}U(t,s)F(s,x)\widetilde{N}(ds,dx)\nonumber\\
&\quad+\int_{a}^{t}\int_{|x|_{V}\geq 1}U(t,s)G(s,x)N(ds,dx).
\end{align*}
Let $\omega(t)=Y(t)-Z(t)$. Then $\omega(t)=U(t,a)\omega(a)$. From the exponential dissipation property of $U(t,s)$, it follows that
$\mathbb{E}\|\omega(t)\|^{2}\leq M^{2}e^{-2\delta(t-a)}\mathbb{E}\|\omega(a)\|^{2}, \ \ t\geq a.$ So $\mathbb{E}\|\omega(t)\|^{2}\rightarrow 0$ as $t\rightarrow \infty$. Since $\omega(t)$ is square-mean almost automorphic, for any sequence of real numbers $\{s'_{n}\}$, there exists a subsequence $\{s_{n}\}$ such that for some stochastic process $\widetilde{\omega}(t)$,
\begin{align}\label{w1}
\lim_{n\rightarrow \infty}\mathbb{E}\|\omega(t+s_{n})-\widetilde{\omega}(t)\|^{2}=0,
\end{align}
and
\begin{align}\label{w2}
\lim_{n\rightarrow \infty}\mathbb{E}\|\widetilde{\omega}(t-s_{n})-\omega(t)\|^{2}=0
\end{align}
for each $t\in \mathbb{R}$. In particular, if $\lim_{n\rightarrow \infty}s'_{n}=\infty$, by (\ref{w1}), we have $\mathbb{E}\|\widetilde{\omega}(t)\|^{2}=0$ for each $t\in \mathbb{R}$. Hence by (\ref{w2}), we have $\mathbb{E}\|{\omega}(t)\|^{2}=0$ for each $t\in \mathbb{R}$. Therefore, $Y(a)=Z(a)$ almost surely, a contradiction. The proof is complete.  \ \ \ $\Box$
\section{Square-mean weighted pseudo almost automorphic solutions of nonautonomous linear stochastic differential equations}
In this section, we establish the existence and uniqueness of square-mean weighted pseudo almost automorphic solutions for nonautonomous linear stochastic differential equations (\ref{levy}). \\
\textbf{Theorem 4.1.} If the assumptions $(H_{1})-(H_{4})$ hold, and $f,g\in SWPAA(\mathbb{R}, \rho)$,  $F,G\in PSWPAA(\mathbb{R}\times V, \rho)$, $\rho\in \mathcal{U}^{inv}\cap \mathcal{U}_{p}^{inv}$,  then the equation (\ref{levy}) has a unique square-mean weighted pseudo almost automorphic mild solution. \\
\textbf{Proof.} It is easy to see that the process (\ref{process}) is a mild solution of (\ref{levy}). Since $f,g\in SWPAA(\mathbb{R}, \rho)$,  $F,G\in PSWPAA(\mathbb{R}\times V, \rho)$, there exist $l\in SAA(\mathbb{R}, \mathcal{L}^{2}(P,H))$, $p\in SAA(\mathbb{R}, L(V,\mathcal{L}^{2}(P,H)))$, $h,\alpha\in PSAA(\mathbb{R}\times V, \mathcal{L}^{2}(P,H)))$ and $m\in SBC_{0}(\mathbb{R}, \rho)$, $q\in SBC_{0}(\mathbb{R}\times V, \rho)$, $\varphi,\beta\in PSBC_{0}(\mathbb{R}\times V, \rho)$ such that
$f(t)=l(t)+m(t)$, $g(t)=p(t)+q(t)$, $F(t,x)=h(t,x)+\varphi(t,x)$,  $G(t,x)=\alpha(t,x)+\beta(t,x)$. We have
\begin{align*}
Y(t)&=\int_{-\infty}^{t}U(t,s)f(s)ds+\int_{-\infty}^{t}U(t,s)g(s)dW(s)\nonumber\\
&\quad+\int_{-\infty}^{t}\int_{|x|_{V}<1}U(t,s)F(s,x)\widetilde{N}(ds,dx)\nonumber\\
&\quad+\int_{-\infty}^{t}\int_{|x|_{V}\geq1}U(t,s)G(s,x)N(ds,dx)\nonumber\\
&=[\int_{-\infty}^{t}U(t,s)l(s)ds+\int_{-\infty}^{t}U(t,s)p(s)dW(s)\nonumber\\
&\quad+\int_{-\infty}^{t}\int_{|x|_{V}<1}U(t,s)h(s,x)\widetilde{N}(ds,dx)\nonumber\\
&\quad+\int_{-\infty}^{t}\int_{|x|_{V}\geq1}U(t,s)\alpha(s,x)N(ds,dx)]\nonumber\\
&\quad+[\int_{-\infty}^{t}U(t,s)m(s)ds+\int_{-\infty}^{t}U(t,s)q(s)dW(s)\nonumber\\
&\quad+\int_{-\infty}^{t}\int_{|x|_{V}<1}U(t,s)\varphi(s,x)\widetilde{N}(ds,dx)\nonumber\\
&\quad+\int_{-\infty}^{t}\int_{|x|_{V}\geq1}U(t,s)\beta(s,x)N(ds,dx)]\nonumber\\
&=:Y_{1}(t)+Y_{2}(t).
\end{align*}
From the proof of Theorem 3.2, we have $Y_{1}\in SAA(\mathbb{R}, \mathcal{L}^{2}(P,H))$.    We need to prove that  $Y_{2}\in SBC_{0}(\mathbb{R},\rho)$.

First, we show that $Y_{2}(t)$ is $\mathcal{L}^{2}$-continuous and $\mathcal{L}^{2}$-bounded. Set $R(t)=\int_{-\infty}^{t}U(t,s)m(s)ds+\int_{-\infty}^{t}U(t,s)q(s)dW(s)$. Similar to
the proof of \cite[Theorem 4.1]{CZ} with minor modifications, it follows that $R(t)$ is $\mathcal{L}^{2}$-continuous and $\mathcal{L}^{2}$-bounded. Since $\varphi,\beta\in PSBC_{0}(\mathbb{R}\times V, \rho)$, then $\varphi,\beta$ are both Poisson stochastically bounded and Poisson stochastically continuous. Similar to the argument for step 1 of the proof of Theorem 3.2, it follows that
\begin{align*}
\int_{-\infty}^{t}\int_{|x|_{V}<1}U(t,s)\varphi(s,x)\widetilde{N}(ds,dx)+\int_{-\infty}^{t}\int_{|x|_{V}\geq1}U(t,s)\beta(s,x)N(ds,dx)
\end{align*}
is $\mathcal{L}^{2}$-continuous.
By properties of the integral for the Poisson random measure, we have
\begin{align}\label{sbc}
&\mathbb{E}\big\|\int_{-\infty}^{t}\int_{|x|_{V}<1}U(t,s)\varphi(s,x)\widetilde{N}(ds,dx)+\int_{-\infty}^{t}\int_{|x|_{V}\geq1}U(t,s)\beta(s,x)N(ds,dx)\big\|^{2}\nonumber\\
&\leq 2E\big\|\int_{-\infty}^{t}\int_{|x|_{V}<1}U(t,s)\varphi(s,x)\widetilde{N}(ds,dx)\big\|^{2}+2\mathbb{E}\big\|\int_{-\infty}^{t}\int_{|x|_{V}\geq1}U(t,s)\beta(s,x)N(ds,dx)\big\|^{2}\nonumber\\
&\leq  2E\big\|\int_{-\infty}^{t}\int_{|x|_{V}<1}U(t,s)\varphi(s,x)\widetilde{N}(ds,dx)\big\|^{2}+4\mathbb{E}\big\|\int_{-\infty}^{t}\int_{|x|_{V}\geq1}U(t,s)\beta(s,x)\widetilde{N}(ds,dx)\big\|^{2}\nonumber\\
&\quad+4\mathbb{E}\big\|\int_{-\infty}^{t}\int_{|x|_{V}\geq1}U(t,s)\beta(s,x)\nu(dx)ds\big\|^{2}\nonumber\\
&\leq  2\int_{-\infty}^{t}\int_{|x|_{V}<1}\|U(t,s)\|^{2}\mathbb{E}\|\varphi(s,x)\|^{2}\nu(dx)ds+4\int_{-\infty}^{t}\int_{|x|_{V}\geq1}\|U(t,s)\|^{2}\mathbb{E}\|\beta(s,x)\|^{2}\nu(dx)ds\nonumber\\  &\quad+4\mathbb{E}\big\|\int_{-\infty}^{t}\int_{|x|_{V}\geq1}U(t,s)\beta(s,x)\nu(dx)ds\big\|^{2}.
\end{align}
By the exponential dissipation property of $U(t,s)$ and $\varphi,\beta\in PSBC_{0}(\mathbb{R}\times V, \mathcal{L}^{2}(P,H))$,  it follows that
\begin{align}\label{sbc1}
&\int_{-\infty}^{t}\int_{|x|_{V}<1}\|U(t,s)\|^{2}\mathbb{E}\|\varphi(s,x)\|^{2}\nu(dx)ds\nonumber\\
&\leq \int_{-\infty}^{t}\int_{|x|_{V}<1}M^{2}e^{-2\delta(t-s)}\mathbb{E}\|\varphi(s,x)\|^{2}\nu(dx)ds\nonumber\\
&\leq M^{2}\int_{-\infty}^{t}e^{-2\delta(t-s)}ds\int_{|x|_{V}<1}\mathbb{E}\|\varphi(s,x)\|^{2}\nu(dx)\nonumber\\
&\leq \frac{M^{2}}{2\delta}\int_{|x|_{V}<1}\mathbb{E}\|\varphi(s,x)\|^{2}\nu(dx)\nonumber\\
&<\infty.
\end{align}
Similarly,
\begin{align}\label{sbc2}
\int_{-\infty}^{t}\int_{|x|_{V}\geq1}\|U(t,s)\|^{2}\mathbb{E}\|\beta(s,x)\|^{2}\nu(dx)ds<\infty.
\end{align}
By the Cauchy-Schwarz inequality, we have
\begin{align}\label{sbc3}
&\mathbb{E}\big\|\int_{-\infty}^{t}\int_{|x|_{V}\geq1}U(t,s)\beta(s,x)\nu(dx)ds\big\|^{2}\nonumber\\
&\leq M^{2}E\left(\int_{-\infty}^{t}\int_{|x|_{V}\geq1}e^{-\frac{\delta(t-s)}{2}}\cdot e^{-\frac{\delta(t-s)}{2}} \cdot \|\beta(s,x)\|\nu(dx)ds\right)^{2}\nonumber\\
&\leq M^{2}\int_{-\infty}^{t}e^{-\delta(t-s)}ds\int_{|x|_{V}\geq 1}\nu(dx)\int_{-\infty}^{t}\int_{|x|_{V}\geq1}e^{-\delta(t-s)}\mathbb{E}\|\beta(s,x)\|^{2}\nu(dx)ds\nonumber\\
&\leq \frac{M^{2}c}{\delta}\int_{-\infty}^{t}e^{-\delta(t-s)}ds\int_{|x|_{V}\geq1}\mathbb{E}\|\beta(s,x)\|^{2}\nu(dx)\nonumber\\
&<\infty.
\end{align}
By (\ref{sbc}), (\ref{sbc1}), (\ref{sbc2}) and (\ref{sbc3}),
\begin{align*}
\mathbb{E}\big\|\int_{-\infty}^{t}\int_{|x|_{V}<1}U(t,s)\varphi(s,x)\widetilde{N}(ds,dx)+\int_{-\infty}^{t}\int_{|x|_{V}\geq1}U(t,s)\beta(s,x)N(ds,dx)\big\|^{2}<\infty.
\end{align*}
Then, $Y_{2}(t)$ is stochastic bounded.

Second, we show that
\begin{align*}
\lim_{r\rightarrow +\infty}\frac{1}{m(r,\rho)}\int_{-r}^{r}\mathbb{E}\|Y_{2}(t)\|^{2}\rho(t)dt=0.
\end{align*}
By the definition of $Y_{2}(t)$, we obtain
\begin{align}\label{y2}
&\frac{1}{m(r,\rho)}\int_{-r}^{r}\mathbb{E}\|Y_{2}(t)\|^{2}\rho(t)dt\nonumber\\
&=\frac{1}{m(r,\rho)}\int_{-r}^{r}\mathbb{E}\big\|\int_{-\infty}^{t}U(t,s)m(s)ds+\int_{-\infty}^{t}U(t,s)q(s)dW(s)\nonumber\\
&\quad+\int_{-\infty}^{t}\int_{|x|_{V}<1}U(t,s)\varphi(s,x)\widetilde{N}(ds,dx)\nonumber\\
&\quad+\int_{-\infty}^{t}\int_{|x|_{V}\geq1}U(t,s)\beta(s,x)N(ds,dx)\big\|^{2}\rho(t)dt\nonumber\\
&\leq \frac{2}{m(r,\rho)}\int_{-r}^{r}\mathbb{E}\big\|\int_{-\infty}^{t}U(t,s)m(s)ds+\int_{-\infty}^{t}U(t,s)q(s)dW(s)\big\|^{2}\rho(t)dt\nonumber\\
&\quad+\frac{2}{m(r,\rho)}\int_{-r}^{r}\mathbb{E}\big\|\int_{-\infty}^{t}\int_{|x|_{V}<1}U(t,s)\varphi(s,x)\widetilde{N}(ds,dx)\nonumber\\
&\quad+\int_{-\infty}^{t}\int_{|x|_{V}\geq1}U(t,s)\beta(s,x)N(ds,dx)\big\|^{2}\rho(t)dt.
\end{align}
Similar to
the proof of \cite[Theorem 4.1]{CZ} with minor modifications, we have
\begin{align}\label{SBC1}
\lim_{r\rightarrow +\infty}\frac{1}{m(r,\rho)}\int_{-r}^{r}\mathbb{E}\big\|\int_{-\infty}^{t}U(t,s)m(s)ds+\int_{-\infty}^{t}U(t,s)q(s)dW(s)\big\|^{2}\rho(t)dt=0.
\end{align}
On the other hand,
\begin{align}\label{SBC2}
&\frac{1}{m(r,\rho)}\int_{-r}^{r}\mathbb{E}\big\|\int_{-\infty}^{t}\int_{|x|_{V}<1}U(t,s)\varphi(s,x)\widetilde{N}(ds,dx)\nonumber\\
&\quad+\int_{-\infty}^{t}\int_{|x|_{V}\geq1}U(t,s)\beta(s,x)N(ds,dx)\big\|^{2}\rho(t)dt\nonumber\\
&\leq \frac{2}{m(r,\rho)}\int_{-r}^{r}\mathbb{E}\big\|\int_{-\infty}^{t}\int_{|x|_{V}<1}U(t,s)\varphi(s,x)\widetilde{N}(ds,dx)\big\|^{2}\rho(t)dt\nonumber\\
&\quad+\frac{2}{m(r,\rho)}\int_{-r}^{r}\mathbb{E}\big\|\int_{-\infty}^{t}\int_{|x|_{V}\geq1}U(t,s)\beta(s,x)N(ds,dx)\big\|^{2}\rho(t)dt\nonumber\\
&\leq \frac{2}{m(r,\rho)}\int_{-r}^{r}\mathbb{E}\big\|\int_{-\infty}^{t}\int_{|x|_{V}<1}U(t,s)\varphi(s,x)\widetilde{N}(ds,dx)\big\|^{2}\rho(t)dt\nonumber\\
&\quad+\frac{4}{m(r,\rho)}\int_{-r}^{r}\mathbb{E}\big\|\int_{-\infty}^{t}\int_{|x|_{V}\geq1}U(t,s)\beta(s,x)\widetilde{N}(ds,dx)\big\|^{2}\rho(t)dt\nonumber\\
&\quad+\frac{4}{m(r,\rho)}\int_{-r}^{r}\mathbb{E}\big\|\int_{-\infty}^{t}\int_{|x|_{V}\geq1}U(t,s)\beta(s,x)\nu(dx)ds\big\|^{2}\rho(t)dt.
\end{align}
By properties of the integral for the Poisson random measure and the exponential dissipation property of $U(t,s)$, we obtain
\begin{align}\label{SBC3}
&\frac{1}{m(r,\rho)}\int_{-r}^{r}\mathbb{E}\big\|\int_{-\infty}^{t}\int_{|x|_{V}<1}U(t,s)\varphi(s,x)\widetilde{N}(ds,dx)\big\|^{2}\rho(t)dt\nonumber\\
&\leq \frac{1}{m(r,\rho)}\int_{-r}^{r}\rho(t)dt\int_{-\infty}^{t}\int_{|x|_{V}<1}\|U(t,s)\|^{2}\mathbb{E}\|\varphi(s,x)\|^{2}\nu(dx)ds\nonumber\\
&\leq \frac{M^{2}}{m(r,\rho)}\int_{-r}^{r}\rho(t)dt\int_{-\infty}^{t}\int_{|x|_{V}<1}e^{-2\delta(t-s)}\mathbb{E}\|\varphi(s,x)\|^{2}\nu(dx)ds\nonumber\\
&=\frac{M^{2}}{m(r,\rho)}\int_{-r}^{r}\rho(t)dt\int_{0}^{+\infty}\int_{|x|_{V}<1}e^{-2\delta s}\mathbb{E}\|\varphi(t-s,x)\|^{2}\nu(dx)ds.
\end{align}
By the Fubini theorem, we get
\begin{align}\label{Fubini}
&\frac{M^{2}}{m(r,\rho)}\int_{-r}^{r}\rho(t)dt\int_{0}^{+\infty}\int_{|x|_{V}<1}e^{-2\delta s}\mathbb{E}\|\varphi(t-s,x)\|^{2}\nu(dx)ds\nonumber\\
&=M^{2}\int_{0}^{+\infty}e^{-2\delta s}ds\frac{1}{m(r,\rho)}\int_{-r}^{r}\int_{|x|_{V}<1}\mathbb{E}\|\varphi(t-s,x)\|^{2}\rho(t)\nu(dx)dt.
\end{align}
From $\rho\in \mathcal{U}^{inv}$, $\varphi\in PSBC_{0}(\mathbb{R}\times V, \rho)$, it follows that
\begin{align}\label{SBC4}
\lim_{r\rightarrow +\infty}\frac{1}{m(r,\rho)}\int_{-r}^{r}\int_{|x|_{V}<1}\mathbb{E}\|\varphi(t-s,x)\|^{2}\rho(t)\nu(dx)dt=0.
\end{align}
By (\ref{SBC3}), (\ref{Fubini}), (\ref{SBC4}) and the Lebesgue dominated convergence theorem, we have
\begin{align}\label{SBC5}
\lim_{r\rightarrow +\infty}\frac{1}{m(r,\rho)}\int_{-r}^{r}\mathbb{E}\big\|\int_{-\infty}^{t}\int_{|x|_{V}<1}U(t,s)\varphi(s,x)\widetilde{N}(ds,dx)\big\|^{2}\rho(t)dt=0.
\end{align}
Similarly,
\begin{align}\label{SBC6}
\lim_{r\rightarrow +\infty}\frac{1}{m(r,\rho)}\int_{-r}^{r}\mathbb{E}\big\|\int_{-\infty}^{t}\int_{|x|_{V}\geq1}U(t,s)\beta(s,x)\widetilde{N}(ds,dx)\big\|^{2}\rho(t)dt=0.
\end{align}
By the Cauchy-Schwarz inequality, the exponential dissipation property of $U(t,s)$ and the Fubini theorem, we have
\begin{align}\label{SBC7}
&\frac{1}{m(r,\rho)}\int_{-r}^{r}\mathbb{E}\big\|\int_{-\infty}^{t}\int_{|x|_{V}\geq1}U(t,s)\beta(s,x)\nu(dx)ds\big\|^{2}\rho(t)dt\nonumber\\
&=\frac{1}{m(r,\rho)}\int_{-r}^{r}\rho(t)dt\ \mathbb{E}\big\|\int_{-\infty}^{t}\int_{|x|_{V}\geq1}U(t,s)\beta(s,x)\nu(dx)ds\big\|^{2}\nonumber\\
&\leq \frac{M^{2}}{m(r,\rho)}\int_{-r}^{r}\rho(t)dt \ \mathbb{E}\left(\int_{-\infty}^{t}\int_{|x|_{V}\geq1}e^{-\frac{\delta(t-s)}{2}}\cdot e^{-\frac{\delta(t-s)}{2}}\cdot\|\beta(s,x)\|\nu(dx)ds\right)^{2}\nonumber\\
&\leq \frac{M^{2}}{m(r,\rho)}\int_{-r}^{r}\rho(t)dt\left(\int_{-\infty}^{t}e^{-\delta(t-s)}ds\int_{|x|_{V}\geq 1}\nu(dx)\int_{-\infty}^{t}\int_{|x|_{V}\geq1}e^{-\delta(t-s)}\mathbb{E}\|\beta(s,x)\|^{2}\nu(dx)ds\right)\nonumber\\
&\leq \frac{M^{2}c}{\delta m(r,\rho)}\int_{-r}^{r}\rho(t)dt\int_{-\infty}^{t}\int_{|x|_{V}\geq1}e^{-\delta(t-s)}\mathbb{E}\|\beta(s,x)\|^{2}\nu(dx)ds\nonumber\\
&=\frac{M^{2}c}{\delta m(r,\rho)}\int_{-r}^{r}\rho(t)dt \int_{0}^{+\infty}\int_{|x|_{V}\geq1}e^{-\delta s}\mathbb{E}\|\beta(t-s,x)\|^{2}\nu(dx)ds\nonumber\\
&=\frac{M^{2}c}{\delta}\int_{0}^{+\infty}e^{-\delta s}ds\frac{1}{m(r,\rho)}\int_{-r}^{r}\int_{|x|_{V}\geq1}\mathbb{E}\|\beta(t-s,x)\|^{2}\rho(t)\nu(dx)dt
\end{align}
Since $\rho\in \mathcal{U}_{p}^{inv}$, $\beta\in PSBC_{0}(\mathbb{R}\times V, \rho)$, by (\ref{SBC7}) and the Lebesgue dominated convergence theorem,
we obtain
\begin{align}\label{SBC8}
\lim_{r\rightarrow +\infty}\frac{1}{m(r,\rho)}\int_{-r}^{r}\mathbb{E}\big\|\int_{-\infty}^{t}\int_{|x|_{V}\geq1}U(t,s)\beta(s,x)\nu(dx)ds\big\|^{2}\rho(t)dt=0.
\end{align}
By (\ref{SBC2}), (\ref{SBC5}), (\ref{SBC6}) and (\ref{SBC8}), we get
\begin{align}\label{SBC9}
&\lim_{r\rightarrow +\infty}\frac{1}{m(r,\rho)}\int_{-r}^{r}\mathbb{E}\big\|\int_{-\infty}^{t}\int_{|x|_{V}<1}U(t,s)\varphi(s,x)\widetilde{N}(ds,dx)\nonumber\\
&\quad+\int_{-\infty}^{t}\int_{|x|_{V}\geq1}U(t,s)\beta(s,x)N(ds,dx)\big\|^{2}\rho(t)dt=0
\end{align}
By (\ref{y2}), (\ref{SBC1}) and (\ref{SBC9}), we have
\begin{align*}
\lim_{r\rightarrow +\infty}\frac{1}{m(r,\rho)}\int_{-r}^{r}\mathbb{E}\|Y_{2}(t)\|^{2}\rho(t)dt=0.
\end{align*}
Therefore, $Y_{2}\in SBC_{0}(\mathbb{R},\rho)$. The proof is complete.  \ \ $\Box$ \\
\section{Square-mean weighted pseudo almost automorphic solutions of nonautonomous semilinear stochastic differential equations}
In this section, we consider the following nonautonomous semilinear stochastic differential equation
\begin{align}\label{nonautonmous}
dY(t)&=A(t)Y(t-)dt+f(t,Y(t-))dt+g(t,Y(t-))dW(t)\nonumber\\
&\quad+\int_{|x|_{V}<1}F(t,Y(t-),x)\widetilde{N}(dt,dx)+\int_{|x|_{V}\geq 1}G(t,Y(t-),x)N(dt,dx),
\end{align}
where $f: \mathbb{R}\times \mathcal{L}^{2}(P,H)\rightarrow \mathcal{L}^{2}(P,H)$, $g: \mathbb{R}\times \mathcal{L}^{2}(P,H)\rightarrow L(V, \mathcal{L}^{2}(P,H))$, $F,G: \mathbb{R}\times \mathcal{L}^{2}(P,H)\times V\rightarrow \mathcal{L}^{2}(P,H)$, $W$ and $N$ are the same as in the previous section, $A(t)$ satisfies assumptions $(H_{1})-(H_{4})$. \\
\textbf{Definition 5.1.} An $\mathcal{F}_{t}$-progressively measurable stochastic process $\{Y(t)\}_{t\in \mathbb{R}}$ is called a mild solution of (\ref{nonautonmous}) if it satisfies the corresponding stochastic integral equation
\begin{align*}
Y(t)&=U(t,a)Y(a)+\int_{a}^{t}U(t,s)f(s,Y(s-))ds+\int_{a}^{t}U(t,s)g(s,Y(s-))dW(s)\nonumber\\
&\quad+\int_{a}^{t}\int_{|x|_{V}<1}U(t,s)F(s,Y(s-),x)\widetilde{N}(ds,dx)\nonumber\\
&\quad+\int_{a}^{t}\int_{|x|_{V}\geq 1}U(t,s)G(s,Y(s-),x)N(ds,dx),
\end{align*}
for all $t\geq a$ and each $a\in \mathbb{R}$.\\
\textbf{Theorem 5.2.} Assume $f,g\in SWPAA(\mathbb{R}\times \mathcal{L}^{2}(P,H),\rho)$, $F=\varphi_{1}+\phi_{1}\in PSWPAA(\mathbb{R}\times \mathcal{L}^{2}(P,H)\times V, \rho)$ with $\varphi_{1}\in PSAA(\mathbb{R}\times \mathcal{L}^{2}(P,H)\times V,\mathcal{L}^{2}(P,H))$ and $\phi_{1}\in PSBC_{0}(\mathbb{R}\times \mathcal{L}^{2}(P,H)\times V, \rho)$, $G=\varphi_{2}+\phi_{2}\in PSWPAA(\mathbb{R}\times \mathcal{L}^{2}(P,H)\times V, \rho)$ with $\varphi_{2}\in PSAA(\mathbb{R}\times \mathcal{L}^{2}(P,H)\times V,\mathcal{L}^{2}(P,H))$ and $\phi_{2}\in PSBC_{0}(\mathbb{R}\times \mathcal{L}^{2}(P,H)\times V, \rho)$, where $\rho\in \mathcal{U}^{inv}\cap \mathcal{U}_{p}^{inv}$. In addition, suppose that $f$, $g$, $F$ and $G$ satisfy the Lipschitz conditions in $Y$ uniformly for $t$, that is, for all $Y,Z\in \mathcal{L}^{2}(P,H)$ and $t\in \mathbb{R}$,
\begin{align*}
&\mathbb{E}\|f(t,Y)-f(t,Z)\|^{2}\leq L\mathbb{E}\|Y-Z\|^{2}, \\
&\mathbb{E}\|(g(t,Y)-g(t,Z))Q^{1/2}\|^{2}_{L(V, \mathcal{L}^{2}(P,H))}\leq L\mathbb{E}\|Y-Z\|^{2}, \\
&\int_{|x|_{V}<1}\mathbb{E}\|F(t,Y,x)-F(t,Z,x)\|^{2}\nu(dx)\leq L\mathbb{E}\|Y-Z\|^{2}, \\
&\int_{|x|_{V}<1}\mathbb{E}\|\varphi_{1}(t,Y,x)-\varphi_{1}(t,Z,x)\|^{2}\nu(dx)\leq L\mathbb{E}\|Y-Z\|^{2}, \\
&\int_{|x|_{V}\geq 1}\mathbb{E}\|G(t,Y,x)-G(t,Z,x)\|^{2}\nu(dx)\leq L\mathbb{E}\|Y-Z\|^{2},\\
&\int_{|x|_{V}\geq 1}\mathbb{E}\|\varphi_{2}(t,Y,x)-\varphi_{2}(t,Z,x)\|^{2}\nu(dx)\leq L\mathbb{E}\|Y-Z\|^{2},
\end{align*}
for some constant $L>0$ is independent of $t$, where $\frac{1+2c}{\delta^{2}}+\frac{2}{\delta}<\frac{1}{4M^{2}L}$. Then (\ref{nonautonmous}) has a unique square-mean weighted pseudo almost automorphic mild solution.\\
\textbf{Proof.} It is easy to prove that $Y(t)$ is a mild solution of (\ref{nonautonmous}) if and only if it satisfies the following integral equation \begin{align*}
Y(t)&=\int_{-\infty}^{t}U(t,s)f(s,Y(s-))ds+\int_{-\infty}^{t}U(t,s)g(s,Y(s-))dW(s)\\
&\quad+\int_{-\infty}^{t}\int_{|x|_{V}<1}U(t,s)F(s,Y(s-),x)\widetilde{N}(ds,dx)\\
&\quad+\int_{-\infty}^{t}\int_{|x|_{V}\geq 1}U(t,s)G(s,Y(s-),x)N(ds,dx).
\end{align*}
For any $Y\in SWPAA(\mathbb{R},\rho)$, define the nonlinear operator $\mathcal{S}$ by
\begin{align*}
(\mathcal{S}Y)(t):&=\int_{-\infty}^{t}U(t,s)f(s,Y(s-))ds+\int_{-\infty}^{t}U(t,s)g(s,Y(s-))dW(s)\\
&\quad+\int_{-\infty}^{t}\int_{|x|_{V}<1}U(t,s)F(s,Y(s-),x)\widetilde{N}(ds,dx)\\
&\quad+\int_{-\infty}^{t}\int_{|x|_{V}\geq 1}U(t,s)G(s,Y(s-),x)N(ds,dx).
\end{align*}
By Lemma 2.17, $f(t,Y(t))\in SWPAA(\mathbb{R},\rho)$ and $g(t,Y(t))\in SWPAA(\mathbb{R},\rho)$ if $Y\in SWPAA (\mathbb{R}, \rho)$. By Theorem 2.33, we see that $J_{1}(t,x):=F(t,Y(t),x)\in PSWPAA(\mathbb{R}\times V, \rho)$ and $J_{2}(t,x):=G(t,Y(t),x)\in PSWPAA(\mathbb{R}\times V, \rho)$ if $Y\in SWPAA (\mathbb{R}, \rho)$. By Theorem 4.1, we have $(\mathcal{S}Y)(t)\in SWPAA(\mathbb{R},\rho)$. So, $\mathcal{S}$ maps $SWPAA(\mathbb{R},\rho)$ into itself.

Next, we show that $\mathcal{S}$ is a contraction mapping on $SWPAA(\mathbb{R},\rho)$. Similar to the proof of \cite[Theorem 4.2]{YZ} with minor modifications, we can obtain the contraction property of $\mathcal{S}$, so the proof is omitted here.

Therefore, $\mathcal{S}$ has a unique fixed point $Y^{\ast}\in SWPAA(\mathbb{R},\rho)$ with $SY^{\ast}=Y^{\ast}$, which is the unique square-mean weighted pseudo almost automorphic mild solution of (\ref{nonautonmous}), the proof is complete. \ \ $\Box$ \\
\textbf{Theorem 5.3.} Assume $f,g\in SWPAA(\mathbb{R}\times \mathcal{L}^{2}(P,H),\rho)$, $F=\varphi_{1}+\phi_{1}\in PSWPAA(\mathbb{R}\times \mathcal{L}^{2}(P,H)\times V, \rho)$ with $\varphi_{1}\in PSAA(\mathbb{R}\times \mathcal{L}^{2}(P,H)\times V,\mathcal{L}^{2}(P,H))$ and $\phi_{1}\in PSBC_{0}(\mathbb{R}\times \mathcal{L}^{2}(P,H)\times V, \rho)$, $G=\varphi_{2}+\phi_{2}\in PSWPAA(\mathbb{R}\times \mathcal{L}^{2}(P,H)\times V, \rho)$ with $\varphi_{2}\in PSAA(\mathbb{R}\times \mathcal{L}^{2}(P,H)\times V,\mathcal{L}^{2}(P,H))$ and $\phi_{2}\in PSBC_{0}(\mathbb{R}\times \mathcal{L}^{2}(P,H)\times V, \rho)$, where $\rho\in \mathcal{U}^{inv}\cap \mathcal{U}_{p}^{inv}$. In addition, suppose that $f(t,0)=g(t,0)=F(t,0,x)=G(t,0,x)=0$, and $f$, $g$, $F$ and $G$ satisfy the local Lipschitz conditions in $Y$ uniformly for $t$, that is, for all $Y,Z\in B_{r}$ and $t\in \mathbb{R}$,
\begin{align*}
&\mathbb{E}\|f(t,Y)-f(t,Z)\|^{2}\leq L_{r}\mathbb{E}\|Y-Z\|^{2}, \\
&\mathbb{E}\|(g(t,Y)-g(t,Z))Q^{1/2}\|^{2}_{L(V, \mathcal{L}^{2}(P,H))}\leq L_{r}\mathbb{E}\|Y-Z\|^{2}, \\
&\int_{|x|_{V}<1}\mathbb{E}\|F(t,Y,x)-F(t,Z,x)\|^{2}\nu(dx)\leq L_{r}\mathbb{E}\|Y-Z\|^{2}, \\
&\int_{|x|_{V}<1}\mathbb{E}\|\varphi_{1}(t,Y,x)-\varphi_{1}(t,Z,x)\|^{2}\nu(dx)\leq L_{r}\mathbb{E}\|Y-Z\|^{2}, \\
&\int_{|x|_{V}\geq 1}\mathbb{E}\|G(t,Y,x)-G(t,Z,x)\|^{2}\nu(dx)\leq L_{r}\mathbb{E}\|Y-Z\|^{2},\\
&\int_{|x|_{V}\geq 1}\mathbb{E}\|\varphi_{2}(t,Y,x)-\varphi_{2}(t,Z,x)\|^{2}\nu(dx)\leq L_{r}\mathbb{E}\|Y-Z\|^{2},
\end{align*}
where $B_{r}=\{Y\in \mathcal{L}^{2}(P,H)|\|Y\|_{2}\leq r\}$,  $L_{r}>0$ is independent of $t$ and $\frac{1+2c}{\delta^{2}}+\frac{2}{\delta}<\frac{1}{4M^{2}L_{r}}$. Then (\ref{nonautonmous}) has a unique square-mean weighted pseudo almost automorphic mild solution.\\
\textbf{Proof.} For $\rho\in \mathcal{U}^{inv}\cap \ \mathcal{U}_{p}^{inv}$, define $SWPAA(\mathbb{R},B_{r})=\{Y| Y\in SWPAA(\mathbb{R},\rho),\ \|Y\|_{L^{2}}\leq r\}$. Since $SWPAA(\mathbb{R},\rho)$ is complete, then $SWPAA(\mathbb{R},B_{r})$ is complete. For any $Y\in SWPAA(\mathbb{R},\rho)$, define the nonlinear operator $\mathcal{S}$ by
\begin{align*}
(\mathcal{S}Y)(t):&=\int_{-\infty}^{t}U(t,s)f(s,Y(s-))ds+\int_{-\infty}^{t}U(t,s)g(s,Y(s-))dW(s)\\
&\quad+\int_{-\infty}^{t}\int_{|x|_{V}<1}U(t,s)F(s,Y(s-),x)\widetilde{N}(ds,dx)\\
&\quad+\int_{-\infty}^{t}\int_{|x|_{V}\geq 1}U(t,s)G(s,Y(s-),x)N(ds,dx).
\end{align*}
Then
\begin{align}\label{S}
\mathbb{E}\|(\mathcal{S}Y)(t)\|^{2}&\leq 4\mathbb{E}\big\|\int_{-\infty}^{t}U(t,s)f(s,Y(s-))ds\big\|^{2}\nonumber\\
&\quad+4\mathbb{E}\big\|\int_{-\infty}^{t}U(t,s)g(s,Y(s-))dW(s)\big\|^{2}\nonumber\\
&\quad+4\mathbb{E}\big\|\int_{-\infty}^{t}\int_{|x|_{V}<1}U(t,s)F(s,Y(s-),x)\widetilde{N}(ds,dx)\big\|^{2}\nonumber\\
&\quad+4\mathbb{E}\big\|\int_{-\infty}^{t}\int_{|x|_{V}\geq 1}U(t,s)F(s,Y(s-),x)N(ds,dx)\big\|^{2}\nonumber\\
&=\mathcal{S}_{1}+\mathcal{S}_{2}+\mathcal{S}_{3}+\mathcal{S}_{4}.
\end{align}
Since $Y(t)$ is $\mathcal{L}^{2}$-continuous, then $Z(t):=\mathbb{E}\|Y(t)\|^{2}$, $t\in \mathbb{R}$, is continuous.
So $Z(t)=Z(t-)$ for all $t\in \mathbb{R}$. \\
For $\mathcal{S}_{1}$, we have
\begin{align}\label{J1}
&4\mathbb{E}\big\|\int_{-\infty}^{t}U(t,s)f(s,Y(s-))ds\big\|^{2}\nonumber\\
&\leq 4M^{2}\mathbb{E}\left(e^{-\frac{\delta(t-s)}{2}}\cdot e^{-\frac{\delta(t-s)}{2}}\cdot \|f(s,Y(s-))\|ds\right)^{2}\nonumber\\
&\leq 4M^{2}\int_{-\infty}^{t}e^{-\delta(t-s)}ds\int_{-\infty}^{t}e^{-\delta(t-s)}\mathbb{E}\|f(s,Y(s-))\|^{2}ds\nonumber\\
&\leq \frac{4M^{2}}{\delta}\int_{-\infty}^{t}e^{-\delta(t-s)}L_{r}\mathbb{E}\|Y(s)\|^{2}ds\nonumber\\
&\leq \frac{4M^{2}L_{r}r^{2}}{\delta^{2}}.
\end{align}
For $\mathcal{S}_{2}$, we have
\begin{align}\label{J2}
&4\mathbb{E}\big\|\int_{-\infty}^{t}U(t,s)g(s,Y(s-))dW(s)\big\|^{2}\nonumber\\
&\leq 4M^{2}\left(\int_{-\infty}^{t}e^{-2\delta(t-s)}\mathbb{E}\|g(s,Y(s-))Q^{1/2}\|^{2}_{L(V, \mathcal{L}^{2}(P,H))}ds\right)\nonumber\\
&\leq 4M^{2}\int_{-\infty}^{t}e^{-2\delta(t-s)}L_{r}\mathbb{E}\|Y(s)\|^{2}ds\nonumber\\
&\leq \frac{2M^{2}L_{r}r^{2}}{\delta}.
\end{align}
For $\mathcal{S}_{3}$, we have
\begin{align}\label{J3}
&4\mathbb{E}\big\|\int_{-\infty}^{t}\int_{|x|_{V}<1}U(t,s)F(s,Y(s-),x)\widetilde{N}(ds,dx)\big\|^{2}\nonumber\\
&\leq 4M^{2}\int_{-\infty}^{t}\int_{|x|_{V}<1}e^{-2\delta(t-s)}\mathbb{E}\|F(s,Y(s-),x)\|^{2}\nu(dx)ds\nonumber\\
&\leq 4M^{2}\int_{-\infty}^{t}e^{-2\delta(t-s)}dsL_{r}r^{2}\nonumber\\
&\leq \frac{2M^{2}L_{r}r^{2}}{\delta}.
\end{align}
For $\mathcal{S}_{4}$, we have
\begin{align}\label{J4}
&4\mathbb{E}\big\|\int_{-\infty}^{t}\int_{|x|_{V}\geq 1}U(t,s)G(s,Y(s-),x)N(ds,dx)\big\|^{2}\nonumber\\
&\leq 8\mathbb{E}\big\|\int_{-\infty}^{t}\int_{|x|_{V}\geq 1}U(t,s)G(s,Y(s-),x)\widetilde{N}(ds,dx)\big\|^{2}\nonumber\\
&\quad+8\mathbb{E}\big\|\int_{-\infty}^{t}\int_{|x|_{V}\geq 1}U(t,s)G(s,Y(s-),x)\nu(dx)ds\big\|^{2}\nonumber\\
&\leq\frac{4M^{2}L_{r}r^{2}}{\delta}+8M^{2}\mathbb{E}\left(\int_{-\infty}^{t}\int_{|x|_{V}\geq 1}e^{-\frac{\delta(t-s)}{2}}\cdot e^{-\frac{\delta(t-s)}{2}}\|G(s,Y(s-),x)\|\nu(dx)ds\right)^{2}\nonumber\\
&\leq\frac{4M^{2}L_{r}r^{2}}{\delta}+8M^{2}\int_{-\infty}^{t}e^{-\delta(t-s)}ds\int_{|x|_{V}\geq 1}\nu(dx)\int_{-\infty}^{t}\int_{|x|_{V}\geq 1}e^{-\delta(t-s)}\mathbb{E}\|G(s,Y(s-),x)\|^{2}\nu(dx)ds\nonumber\\
&\leq\frac{4M^{2}L_{r}r^{2}}{\delta}+\frac{8M^{2}cL_{r}r^{2}}{\delta^{2}}.
\end{align}
From (\ref{S}), (\ref{J1}), (\ref{J2}), (\ref{J3}) and (\ref{J4}), it follows that
\begin{align*}
\mathbb{E}\|(\mathcal{S}Y)(t)\|^{2}\leq \left(\frac{4M^{2}L_{r}+8M^{2}cL_{r}}{\delta^{2}}+\frac{8M^{2}L_{r}}{\delta}\right)r^{2}<r^{2},
\end{align*}
that is $\mathcal{S}Y\in B_{r}$.

Since $Y\in SWPAA (\mathbb{R}, \rho)$, by Theorem 2.33, $J_{1}(t,x):=F(t,Y(t),x)\in PSWPAA(\mathbb{R}\times V, \rho)$ and $J_{2}(t,x):=G(t,Y(t),x)\in PSWPAA(\mathbb{R}\times V, \rho)$. By Theorem 4.1, we have $(\mathcal{S}Y)(t)\in SWPAA(\mathbb{R},\rho)$. Hence $\mathcal{S}Y\in SWPAA(\mathbb{R},B_{r})$. So, $\mathcal{S}$ maps $SWPAA(\mathbb{R},\rho)$ into itself.  The proof of contraction property of $\mathcal{S}$ is is omitted, since it is similar to the proof of \cite[Theorem 4.2]{YZ} with minor modifications. Therefore, $\mathcal{S}$ has a unique fixed point $Y^{\ast}\in SWPAA(\mathbb{R},B_{r})$ with $SY^{\ast}=Y^{\ast}$, which is the unique square-mean weighted pseudo almost automorphic mild solution of (\ref{nonautonmous}), the proof is complete. \ \ $\Box$ \\
\section{Stability of square-mean weighted pseudo almost automorphic solutions}
In this section, we investigate the stability of square-mean weighted pseudo almost automorphic solutions of nonautonomous nonlinear stochastic differential equation (\ref{nonautonmous}).\\
\textbf{Definition 6.1} The unique square-mean weighted pseudo almost automorphic solution $Y^{\ast}(t)$ of (\ref{nonautonmous}) is said to be
stable in the square-mean sense, if for arbitrary $\varepsilon>0$, there exists $\eta>0$ such that
\begin{align*}
\mathbb{E}\|Y(t)-Y^{\ast}(t)\|^{2}<\varepsilon, \ \ t\geq 0,
\end{align*}
whenever $\mathbb{E}\|Y(0)-Y^{\ast}(0)\|^{2}<\eta$, where $Y(t)$ stands for the solution of (\ref{nonautonmous}) with initial value $Y(0)$. The solution $Y^{\ast}(t)$ is said to be exponentially stable in the square-mean sense if it is stable in the square-mean sense and for some $\varepsilon>0$,
\begin{align}\label{stable}
\lim_{t\rightarrow +\infty}e^{\varepsilon t}\mathbb{E}\|Y(t)-Y^{\ast}(t)\|^{2}=0.
\end{align}
If (\ref{stable}) holds for any $Y(0)\in \mathcal{L}^{2}(P,H)$, then $Y^{\ast}(t)$ is said to be globally exponentially stable in the square-mean sense. \\
\textbf{Theorem 6.2} Assume that all the conditions of Theorem 5.2 hold, and also assume
\begin{align*}
\frac{5M^{2}L(1+2c)}{\delta^{2}}+\frac{10M^{2}L}{\delta}<1,
\end{align*}
then the square-mean weighted pseudo almost automorphic mild solution $Y^{\ast}(t)$ of (\ref{nonautonmous}) is globally exponentially  stable in square-mean sense. \\
\textbf{Proof.}  Assume that $Y(t)$ is any solution of (\ref{nonautonmous}) with the initial value $Y(0)$ with an interval of existence $[0,T')$.
Then for any $t\in [0,T')$, \begin{align*}
&\mathbb{E}\big\|Y(t)-Y^{\ast}(t)\big\|^{2}\\
&=\mathbb{E}\big\|U(t,0)[Y(0)-Y^{\ast}(0)]+\int_{0}^{t}U(t,s)[f(s,Y(s-))-f(s,Y^{\ast}(s-))]ds\\
&\quad+\int_{0}^{t}U(t,s)[g(s,Y(s-))-g(s,Y^{\ast}(s-))]dW(s)\\
&\quad+\int_{0}^{t}\int_{|x|_{V}<1}U(t,s)[F(s,Y(s-),x)-F(s,Y^{\ast}(s-),x)]\widetilde{N}(ds,dx)\\
&\quad+\int_{0}^{t}\int_{|x|_{V}\geq 1}U(t,s)[G(s,Y(s-),x)-G(s,Y^{\ast}(s-),x)]N(ds,dx)\big\|^{2}\\
&\leq 5\mathbb{E}\|U(t,0)[Y(0)-Y^{\ast}(0)]\|^{2}+5\mathbb{E}\big\|\int_{0}^{t}U(t,s)[f(s,Y(s-))-f(s,Y^{\ast}(s-))]ds\big\|^{2}\nonumber\\
&\quad+5\mathbb{E}\big\|\int_{0}^{t}U(t,s)[g(s,Y(s-))-g(s,Y^{\ast}(s-))]dW(s)\big\|^{2}\nonumber\\
&\quad+5\mathbb{E}\big\|\int_{0}^{t}\int_{|x|_{V}<1}U(t,s)[F(s,Y(s-),x)-F(s,Y^{\ast}(s-),x)]\widetilde{N}(ds,dx) \big\|^{2}\nonumber\\
&\quad+5\mathbb{E}\big\|\int_{0}^{t}\int_{|x|_{V}\geq 1}U(t,s)[G(s,Y(s-),x)-G(s,Y^{\ast}(s-),x)]N(ds,dx)\big\|^{2}\nonumber\\
&\leq 5M^{2}e^{-2\delta t}\mathbb{E}\|Y(0)-Y^{\ast}(0)\|^{2}+\frac{5M^{2}L}{\delta}\int_{0}^{t}e^{-\delta (t-s)}\mathbb{E}\|Y(s-)-Y^{\ast}(s-)\|^{2}ds\\
&\quad+5M^{2}L\int_{0}^{t}e^{-2\delta (t-s)}\mathbb{E}\|Y(s-)-Y^{\ast}(s-)\|^{2}ds+5M^{2}L\int_{0}^{t}e^{-2\delta(t-s)}\mathbb{E}\|Y(s-)-Y^{\ast}(s-)\|^{2}ds\\
&\quad+10M^{2}L\int_{0}^{t}e^{-2\delta(t-s)}\mathbb{E}\|Y(s-)-Y^{\ast}(s-)\|^{2}ds\nonumber\\
&\quad+\frac{10M^{2}cL}{\delta}\int_{0}^{t}e^{-\delta(t-s)}\mathbb{E}\|Y(s-)-Y^{\ast}(s-)\|^{2}ds.
\end{align*}
Define $Z(t):=\mathbb{E}\|Y(t)-Y^{\ast}(t)\|^{2}$, $t\geq 0$. Since $Y(t)$, $Y^{\ast}(t)$ are $\mathcal{L}^{2}$-continuous, then $Z(t)$ is continuous.
So $Z(t)=Z(t-)$ for all $t\geq 0$. Then we have
\begin{align*}
&\mathbb{E}\big\|Y(t)-Y^{\ast}(t)\big\|^{2}\\
&\leq 5M^{2}e^{-2\delta t}\mathbb{E}\|Y(0)-Y^{\ast}(0)\|^{2}+\frac{5M^{2}L}{\delta}\int_{0}^{t}e^{-\delta (t-s)}\mathbb{E}\|Y(s)-Y^{\ast}(s)\|^{2}ds\\
&\quad+5M^{2}L\int_{0}^{t}e^{-2\delta (t-s)}\mathbb{E}\|Y(s)-Y^{\ast}(s)\|^{2}ds+5M^{2}L\int_{0}^{t}e^{-2\delta(t-s)}\mathbb{E}\|Y(s)-Y^{\ast}(s)\|^{2}ds\\
&\quad+10M^{2}L\int_{0}^{t}e^{-2\delta(t-s)}\mathbb{E}\|Y(s)-Y^{\ast}(s)\|^{2}ds+\frac{10M^{2}cL}{\delta}\int_{0}^{t}e^{-\delta(t-s)}\mathbb{E}\|Y(s)-Y^{\ast}(s)\|^{2}ds\\
&\leq 5M^{2}e^{-2\delta t}\mathbb{E}\|Y(0)-Y^{\ast}(0)\|^{2}+\frac{5M^{2}L(1+2c)}{\delta}\int_{0}^{t}e^{-\delta(t-s)}\mathbb{E}\|Y(s)-Y^{\ast}(s)\|^{2}ds\\
&\quad+20M^{2}L\int_{0}^{t}e^{-2\delta(t-s)}\mathbb{E}\|Y(s)-Y^{\ast}(s)\|^{2}ds.
\end{align*}
Since $\frac{5M^{2}L(1+2c)}{\delta^{2}}+\frac{10M^{2}L}{\delta}<1$ as assumed above, there exists a number $\varepsilon\in (0,\delta)$ such that
\begin{align*}
\frac{5M^{2}L(1+2c)}{\delta(\delta-\varepsilon)}+\frac{20M^{2}L}{2\delta-\varepsilon}<1.
\end{align*}
Denote $\mu=\frac{5M^{2}L(1+2c)}{\delta(\delta-\varepsilon)}+\frac{20M^{2}L}{2\delta-\varepsilon}$.
For any $T\in [0,T')$,
\begin{align}\label{estimate}
&\int_{0}^{T}e^{\varepsilon t}\mathbb{E}\|Y(t)-Y^{\ast}(t)\|^{2}dt\nonumber\\
&\leq 5M^{2}\mathbb{E}\|Y(0)-Y^{\ast}(0)\|^{2}\int_{0}^{T}e^{\varepsilon t}e^{-2\delta t}dt\nonumber\\
&\quad+\frac{5M^{2}L(1+2c)}{\delta}\int_{0}^{T}e^{\varepsilon t}dt\int_{0}^{t}e^{-\delta(t-s)}\mathbb{E}\|Y(s)-Y^{\ast}(s)\|^{2}ds\nonumber\\
&\quad+20M^{2}L\int_{0}^{T}e^{\varepsilon t}dt\int_{0}^{t}e^{-2\delta(t-s)}\mathbb{E}\|Y(s)-Y^{\ast}(s)\|^{2}ds\nonumber\\
&\leq \frac{5M^{2}}{2\delta-\varepsilon}\mathbb{E}\|Y(0)-Y^{\ast}(0)\|^{2}\nonumber\\
&\quad+\frac{5M^{2}L(1+2c)}{\delta}\int_{0}^{T}\mathbb{E}\|Y(s)-Y^{\ast}(s)\|^{2}ds\int_{s}^{T}e^{\varepsilon t}e^{-\delta(t-s)}dt\nonumber\\
&\quad+20M^{2}L\int_{0}^{T}\mathbb{E}\|Y(s)-Y^{\ast}(s)\|^{2}ds\int_{s}^{T}e^{\varepsilon t}e^{-2\delta(t-s)}dt\nonumber\\
&\leq \frac{5M^{2}}{2\delta-\varepsilon}\mathbb{E}\|Y(0)-Y^{\ast}(0)\|^{2}\nonumber\\
&\quad+\frac{5M^{2}L(1+2c)}{\delta(\delta-\varepsilon)}\int_{0}^{T}e^{\varepsilon s}\mathbb{E}\|Y(s)-Y^{\ast}(s)\|^{2}ds\nonumber\\
&\quad+\frac{20M^{2}L}{2\delta-\varepsilon}\int_{0}^{T}e^{\varepsilon s}\mathbb{E}\|Y(s)-Y^{\ast}(s)\|^{2}ds.
\end{align}
From (\ref{estimate}), it follows that
\begin{align}\label{exp}
\int_{0}^{T}e^{\varepsilon t}\mathbb{E}\|Y(t)-Y^{\ast}(t)\|^{2}dt\leq \frac{1}{1-\mu}\cdot \frac{5M^{2}}{2\delta-\varepsilon}\mathbb{E}\|Y(0)-Y^{\ast}(0)\|^{2}.
\end{align}
This estimation is independent of $T$ and $T'$, then the maximal interval of the existence of the solution of (\ref{nonautonmous}) can be extended to the positive infinity. By (\ref{exp}), we have
\begin{align*}
\mathbb{E}\|Y(t)-Y^{\ast}(t)\|^{2}\rightarrow 0 \ \mbox{as} \ t\rightarrow +\infty.
\end{align*}
We obtain that the square-mean pseudo almost automorphic mild solution $Y^{\ast}(t)$ of (\ref{nonautonmous}) is globally exponentially stable in the square-mean sense. The proof is complete.  \ $\Box$ \\
\textbf{Theorem 6.3} Assume that all the assumptions of theorem 5.2 hold and that
\begin{align*}
\frac{5M^{2}L_{r}(1+2c)}{\delta^{2}}+\frac{10M^{2}L_{r}}{\delta}<1,
\end{align*}
then the square-mean weighted pseudo almost automorphic mild solution $Y^{\ast}(t)$ of (\ref{nonautonmous}) is locally exponentially stable in the square-mean sense with the attractive domain $D=\{Y \in\mathcal{L}^{2}(P,H)| \|Y\|_{2}\leq r_{1}\}$, where
\begin{align*}
r_{1}< \min\left\{r,\frac{r}{\sqrt{5}M}\sqrt{1-5M^{2}L_{r}\left(\frac{2}{\delta}+\frac{1+2c}{\delta^{2}}\right)}\right\}
\end{align*}
\textbf{Proof.} Assume that $Y(t)$ is any solution of (\ref{nonautonmous}) with the initial value $\|Y(0)\|_{2}\leq r_{1}$. We claim that $\|Y(t)\|_{2}\leq r$ for all $t\geq 0$. If this is not true, then there exists $t_{0}>0$ such that $\|Y(t_{0})\|_{2}=r$ and $\|Y(t)\|_{2}<r$ for $t\in[0,t_{0})$. We have
\begin{align*}
\mathbb{E}\big\|Y(t_{0})\big\|^{2}&=\mathbb{E}\big\|U(t_{0},0)Y(0)+\int_{0}^{t_{0}}U(t_{0},s)f(s,Y(s-))ds+\int_{0}^{t_{0}}U(t_{0},s)g(s,Y(s-))dW(s)\nonumber\\
&\quad+\int_{0}^{t_{0}}\int_{|x|_{V}<1}U(t_{0},s)F(s,Y(s-),x)\widetilde{N}(ds,dx)\nonumber\\
&\quad+\int_{0}^{t_{0}}\int_{|x|_{V}\geq 1}U(t_{0},s)G(s,Y(s-),x)N(ds,dx)\big\|^{2}\nonumber\\
&\leq 5\mathbb{E}\|U(t_{0},0)Y(0)\|^{2}+5\mathbb{E}\big\|\int_{0}^{t_{0}}U(t_{0},s)f(s,Y(s-))ds\big\|^{2}\nonumber\\
&\quad+5\mathbb{E}\big\|\int_{0}^{t_{0}}U(t_{0},s)g(s,Y(s-))dW(s)\big\|^{2}\nonumber\\
&\quad+5\mathbb{E}\big\|\int_{0}^{t_{0}}\int_{|x|_{V}<1}U(t_{0},s)F(s,Y(s-),x)\widetilde{N}(ds,dx) \big\|^{2}\nonumber\\
&\quad+5\mathbb{E}\big\|\int_{0}^{t_{0}}\int_{|x|_{V}\geq 1}U(t_{0},s)G(s,Y(s-),x)N(ds,dx)\big\|^{2}\nonumber\\
&\leq 5M^{2}e^{-2\delta t_{0} }\mathbb{E}\|Y(0)\|^{2}+5M^{2}\int_{0}^{t_{0}}e^{-\delta(t_{0}-s)}ds\int_{0}^{t_{0}}e^{-\delta(t_{0}-s)}\mathbb{E}\|f(s,Y(s-))\|^{2}ds\nonumber\\
&\quad+5M^{2}\left(\int_{0}^{t_{0}}e^{-2\delta(t_{0}-s)}\mathbb{E}\|g(s,Y(s-))Q^{1/2}\|^{2}_{L(V, \mathcal{L}^{2}(P,H))}ds\right)\nonumber\\
&\quad+5M^{2}\int_{0}^{t_{0}}\int_{|x|_{V}<1}e^{-2\delta(t_{0}-s)}\mathbb{E}\|F(s,Y(s-),x)\|^{2}\nu(dx)ds\nonumber\\
&\quad+10\mathbb{E}\big\|\int_{0}^{t_{0}}\int_{|x|_{V}\geq 1}U(t_{0},s)G(s,Y(s-),x)\widetilde{N}(ds,dx)\big\|^{2}\nonumber\\
&\quad+10\mathbb{E}\big\|\int_{0}^{t_{0}}\int_{|x|_{V}\geq 1}U(t_{0},s)G(s,Y(s-),x)\nu(dx)ds\big\|^{2}\nonumber\\
&\leq 5M^{2}r_{1}^{2}+\frac{5M^{2}L_{r}r^{2}}{\delta^{2}}+5M^{2}\int_{0}^{t_{0}}e^{-2\delta(t_{0}-s)}L_{r}\mathbb{E}\|Y(s)\|^{2}ds\nonumber\\
&\quad+5M^{2}\int_{0}^{t_{0}}e^{-2\delta(t_{0}-s)}L_{r}\mathbb{E}\|Y(s)\|^{2}ds+10M^{2}\int_{0}^{t_{0}}e^{-2\delta(t_{0}-s)}L_{r}\mathbb{E}\|Y(s)\|^{2}ds\nonumber\\
&\quad+10M^{2}\mathbb{E}\left(\int_{0}^{t_{0}}\int_{|x|_{V}\geq 1}e^{-\frac{\delta(t_{0}-s)}{2}}\cdot e^{-\frac{\delta(t_{0}-s)}{2}}\|G(s,Y(s-),x)\|\nu(dx)ds\right)^{2}\nonumber\\
&\leq 5M^{2}r_{1}^{2}+\frac{5M^{2}L_{r}r^{2}}{\delta^{2}}+\frac{5M^{2}L_{r}r^{2}}{2\delta}+\frac{5M^{2}L_{r}r^{2}}{2\delta}+\frac{5M^{2}L_{r}r^{2}}{\delta}+\frac{10M^{2}cL_{r}r^{2}}{\delta^{2}}\nonumber\\
&<r^{2},
\end{align*}
this is a contradiction.

By the similar argument as the proof of Theorem 6.2, for $\varepsilon>0$, we have
\begin{align*}
e^{\varepsilon t}\mathbb{E}\|Y(t)-Y^{\ast}(t)\|^{2}dt\rightarrow 0 \ \mbox{as} \ t\rightarrow +\infty.
\end{align*}
Therefore the square-mean pseudo almost automorphic mild solution $Y^{\ast}(t)$ of (\ref{nonautonmous}) is locally exponentially stable in the square-mean sense with the attractive domain $D=\{Y \in\mathcal{L}^{2}(P,H)| \|Y\|_{2}\leq r_{1}\}$. The proof is complete.  \ $\Box$
\section{An example}
Consider the stochastic heat equation with the Dirichlet boundary conditions
\begin{align}\label{example}
\frac{\partial u(t,\xi)}{\partial t}&=\left(\frac{\partial^{2} u}{\partial \xi^{2}}(t,\xi)+u(t,\xi)\sin\frac{1}{2+\sin t+\sin \pi t}\right)+f(t,u(t,\xi))\nonumber\\
&\quad+g(t,u(t,\xi))\frac{\partial W}{\partial t}(t,\xi)+r(t,u(t,\xi))\frac{\partial Z}{\partial t}(t,\xi), \ t\in \mathbb{R},\ \xi\in (0,1), \nonumber\\
&u(t,0)=u(t,1)=0, \ t\in \mathbb{R},
\end{align}
where
\begin{align*}
f(t,u)&=a_{1}u\sin\frac{1}{2+\cos t+\cos \sqrt{2}t}+a_{2}\eta(t)\cos u, \\
g(t,u)&=a_{3}u\sin\frac{1}{2+\cos t+\cos \sqrt{3}t}+a_{4}\eta(t)\sin u, \\
r(t,u)&=a_{5}u\sin\frac{1}{2+\cos t+\cos \sqrt{2}t}+a_{6}\eta(t)\sin u,
\end{align*}
$a_{i} (i=1,\ldots,6)$ are positive constants, $\eta(t)=t\cdot 1_{[0,1]}(t)+1_{[1,\infty)}(t)$, $1_{I}(\cdot)$ is a characteristic function on the interval $I$,  $W$ is a $Q$-Wiener process on $L^{2}[0,1]$ with $\mbox{Tr}Q<\infty$, $Z$ is a L$\acute{e}$vy pure jump process on $L^{2}[0,1]$ which is independent of $W$.
Let $A$ be the Laplace operator, and
\begin{align*}
D(A)=\big\{\tau\in C^{1}[0,1]\big| &\tau'(r) \ \mbox{is absolutely continuous on [0,1]}, \\
&\tau''(r)\in L^{2}[0,1], \ \tau(0)=\tau(1)=0\big\}.
\end{align*}
Then $A$ generates a $C_{0}$-semigroup $(T(t))_{t\geq 0}$ on $L^{2}[0,1]$, which is
\begin{align*}
(T(t)\tau)(r)=\sum_{n=1}^{\infty}e^{-n^{2}\pi^{2}t}(\tau,e_{n})e_{n}(r),
\end{align*}
where $(\cdot,\cdot)$ denotes the inner product on $L^{2}[0,1]$, $e_{n}(r)=\sqrt{2}\sin(n\pi r)$, $n=1,2,\ldots$, and $\|T(t)\|\leq e^{-\pi^{2}t}$ for $t\geq 0$.

Define a family of linear operators $A_{1}(t)$ by
\begin{align*}
&D(A_{1}(t))=D(A), \ t\in \mathbb{R}, \\
&A_{1}(t)\tau=(A+\sin\frac{1}{2+\sin t+ \sin \pi t})\tau, \ \tau\in D(A_{1}(t)).
\end{align*}
Then $(A_{1}(t))_{t\in \mathbb{R}}$ generates an evolution family $\{U(t,s)\}_{t\geq s}$ such that
\begin{align*}
U(t,s)\tau=T(t-s)e^{\int_{s}^{t}\sin \frac{1}{2+\sin \theta+\sin \pi\theta}d\theta}\tau,\ \tau\in \mathcal{L}^{2}(P,L^{2}[0,1]).
\end{align*}
Denote $H=V=L^{2}[0,1]$. The equation (\ref{example}) can be rewritten as the abstract form
\begin{align*}
dY&=A(t)Ydt+F(t,Y)dt+G(t,Y)dW+\int_{|z|_{V}<1}R(t,Y,z)\widetilde{N}(dt,dz)\\
&\quad+\int_{|z|_{V}\geq 1}R(t,Y,z)N(dt,dz)
\end{align*}
on the Hilbert space $H$, where
$F(t,Y):=f(t,u)$, $G(t,Y)dW:=g(t,u)dW$,
\begin{align*}
r(t,u)dZ:=\int_{|z|_{V}<1}R(t,Y,z)\widetilde{N}(dt,dz)+\int_{|z|_{V}\geq 1}R(t,Y,z)N(dt,dz)
\end{align*}
with
\begin{align*}
Z(t,\xi)=\int_{|z|_{V<1}}z\widetilde{N}(t,dz)+\int_{|z|_{V\geq 1}}z N(t,dz),  \ \  R(t,Y,z)=h(t,u)z.
\end{align*}
For simplicity, we assume that the L$\acute{e}$vy pure jump process $Z$ on $L^{2}[0,1]$ is decomposed as above by the L$\acute{e}$vy-It$\hat{o}$ decomposition theorem.

By $\|T(t)\|\leq e^{-\pi^{2}t}$, we have $\|U(t,s)\|\leq e^{-(\pi^{2}-1)(t-s)}$, $t\geq s$. We can choose $M=1$ and $\delta=\pi^{2}-1$. By the almost automorphic property of  $\sin\frac{1}{2+\sin\theta+\sin\pi\theta}$ and
\begin{align*}
U(t+s_{n},s+s_{n})\tau&=T(t-s)e^{\int_{s+s_{n}}^{t+s_{n}}\sin \frac{1}{2+\sin \theta+\sin \pi\theta}d\theta}\tau\\
&=T(t-s)e^{\int_{s}^{t}\sin \frac{1}{2+\sin(\theta+s_{n})+\sin\pi(\theta+s_{n})}d\theta}\tau, \ \tau\in \mathcal{L}^{2}(P, H),
\end{align*}
we see that $U(t,s)\tau\in SBAA(\mathbb{R}\times \mathbb{R}, \mathcal{L}^{2}(P, H))$ uniformly for all $\tau$ in any bounded subset of $\mathcal{L}^{2}(P, H)$. Similarly, $U(t,s)y\in SBAA(\mathbb{R}\times \mathbb{R}, \mathcal{L}^{2}(P,H))$ uniformly for all $y\in \mathbb{B}$, where
\begin{align*}
\mathbb{B}:=\{y: \mathbb{R}\times \mathcal{L}^{2}(P,H) \times V\rightarrow \mathcal{L}^{2}(P,H), \ \sup_{t\in \mathbb{R}}\int_{V}\mathbb{E}\|y(t,Y,x)\|^{2}\nu(dx)<\infty\}.
\end{align*}
Choose $\rho(t)=e^{-t}$, then $\rho\in \mathcal{U}^{inv}\cap \ \mathcal{U}_{p}^{inv}$. It is easy to show that both
\begin{align*}
a_{1}u\sin\frac{1}{2+\cos t+\cos \sqrt{2}t}+a_{2}\eta(t)\cos u
\end{align*}
and
\begin{align*}
a_{3}u\sin\frac{1}{2+\cos t+\cos \sqrt{3}t}+a_{4}\eta(t)\sin u
\end{align*}
are square-mean weighted pseudo almost automprphic processes about $\rho(t)=e^{-t}$, we have that $F\in SWPAA(\mathbb{R}\times \mathcal{L}^{2}(P,H),\rho)$ and $G\in SWPAA(\mathbb{R}\times \mathcal{L}^{2}(P,H),\rho)$. On the other hand, assume that $\nu$ is a finite measure, it is easy to show that  $R\in PSWPAA(\mathbb{R}\times \mathcal{L}^{2}(P,H)\times V,\rho)$. Note that  $F$, $G$, $R$ satisfy the global Lipschitz conditions, with Lipschitz constant $L=\max\{2a_{1}^{2}+2a_{2}^{2}, 2a_{3}^{2}+2a_{4}^{2}, 2a_{5}^{2}+2a_{6}^{2}\}$. If $L$ satisfies the inequality $\frac{5M^{2}L(1+2c)}{\delta^{2}}+\frac{10M^{2}L}{\delta}<1$, then by Theorem 5.2 and Theorem 6.2, the stochastic heat equation (\ref{example}) has a unique square-mean weighted pseudo almost automorphic solution, which is globally exponentially stable in the square-mean sense. \\\\
\textbf{Appendix}\\\\
\textbf{Proof of Lemma 2.32.} (Necessity). For any $\varepsilon>0$,
 \begin{align*}
 &\frac{1}{m(r,\rho)}\int_{-r}^{r}\int_{V}\mathbb{E}\|J(t,x)\|^{2}\nu(dx)\rho(t)dt\\
 &=\frac{1}{m(r,\rho)}\left(\int_{M_{r,\varepsilon}(\varphi)}\int_{V}\mathbb{E}\|J(t,x)\|^{2}\nu(dx)\rho(t)dt+\int_{[-r,r]-M_{r,\varepsilon}(\varphi)}\int_{V}\mathbb{E}\|J(t,x)\|^{2}\nu(dx)\rho(t)dt\right)\\
 &\geq\frac{\varepsilon}{m(r,\rho)}\int_{M_{r,\varepsilon}(\varphi)}\rho(t)dt.
\end{align*}
From $J\in PSBC_{0}(\mathbb{R}\times V, \rho)$, it follows that
\begin{align*}
\lim_{r\rightarrow +\infty}\frac{1}{m(r,\rho)}\int_{-r}^{r}\int_{V}\mathbb{E}\|J(t,x)\|^{2}\nu(dx)\rho(t)dt=0.
\end{align*}
Then
\begin{align*}
\lim_{r\rightarrow +\infty}\frac{1}{m(r,\rho)}\int_{M_{r,\varepsilon}(\varphi)}\rho(t)dt=0.
\end{align*}
(Sufficiency). Since $J\in PSBC_{0}(\mathbb{R}\times V, \rho)$, there exists a constant $M>0$ such that
\begin{align*}
\int_{V}\mathbb{E}\|J(t,x)\|^{2}\nu(dx)\leq M,
\end{align*}
for all $t\in \mathbb{R}$.
From
\begin{equation*}
\lim_{r\rightarrow +\infty}\frac{1}{m(r,\rho)}\int_{M_{r,\varepsilon}(\varphi)}\rho(t)dt=0,
\end{equation*}
we have for any $\varepsilon>0$, there exists a constant $r_{0}$ such that for $r>r_{0}$,
\begin{equation*}
\frac{1}{m(r,\rho)}\int_{M_{r,\varepsilon}(\varphi)}\rho(t)dt<\varepsilon.
\end{equation*}
Hence
\begin{align*}
&\frac{1}{m(r,\rho)}\int_{-r}^{r}\int_{V}\mathbb{E}\|J(t,x)\|^{2}\nu(dx)\rho(t)dt\\
&=\frac{1}{m(r,\rho)}\left(\int_{M_{r,\varepsilon}(\varphi)}\int_{V}\mathbb{E}\|J(t,x)\|^{2}\nu(dx)\rho(t)dt+\int_{[-r,r]-M_{r,\varepsilon}(\varphi)}\int_{V}\mathbb{E}\|J(t,x)\|^{2}\nu(dx)\rho(t)dt\right)\\
&\leq \frac{M}{m(r,\rho)}\int_{M_{r,\varepsilon}(\varphi)}\rho(t)dt+\frac{\varepsilon}{m(r,\rho)}\int_{[-r,r]-M_{r,\varepsilon}(\varphi)}\rho(t)dt\\
&\leq M\varepsilon+\varepsilon.
\end{align*}
Therefore
\begin{align*}
\lim_{r\rightarrow +\infty}\frac{1}{m(r,\rho)}\int_{-r}^{r}\int_{V}\mathbb{E}\|J(t,x)\|^{2}\nu(dx)\rho(t)dt=0.
\end{align*}
This implies that $J\in PSBC_{0}(\mathbb{R}\times V, \rho)$. \ \ $\Box$ \\
\textbf{Proof of Theorem 2.33.} Since $Y\in SWPAA (\mathbb{R}, \rho)$, there exist $\alpha\in SAA(\mathbb{R}, \mathcal{L}^{2}(P,H))$ and $\beta\in SBC_{0}(\mathbb{R},\rho)$ such $Y=\alpha+\beta$. Since $\alpha\in SAA(\mathbb{R}, \mathcal{L}^{2}(P,H))$, and $\int_{V}\mathbb{E}\|g(t,Y,x)-g(t,Z,x)\|^{2}\nu(dx)\leq L\mathbb{E}\|Y-Z\|^{2}$, from Lemma 2.27, it follows that $H\in PSAA(\mathbb{R}\times V,\mathcal{L}^{2}(P,H))$. Note that
  \begin{align*}
  J(t,x)&=g(t,\alpha(t),x)+ J(t,x)-g(t,\alpha(t),x)\\
  &=g(t,\alpha(t),x)+ F(t,Y(t),x)-F(t,\alpha(t),x)+\varphi(t,\alpha(t),x),
  \end{align*}
  let $H(t,x)=g(t,\alpha(t),x)$, $\Phi(t,x)= F(t,Y(t),x)-F(t,\alpha(t),x)+\varphi(t,\alpha(t),x)$, to prove that $J(t,x):=F(t,Y(t),x)$ is Poisson square-mean weighted  pseudo almost automorphic, it is enough to prove $\Phi(t,x)\in PSBC_{0}(\mathbb{R}\times V, \rho)$.

  Firstly, we will prove that $F(t,Y(t),x)-F(t,\alpha(t),x)\in PSBC_{0}(\mathbb{R}\times V, \rho)$. Since $F\in PSWPAA(\mathbb{R}\times \mathcal{L}^{2}(P,H)\times V, \mathcal{L}^{2}(P,H))$, by Definition 2.20 and the $\mathcal{L}^{2}$-continuity of $Y(t)$ and $\alpha(t)$, it is easy to show that $F(t,Y(t),x)-F(t,\alpha(t),x)\in PSC(\mathbb{R}\times V, \mathcal{L}^{2}(P,H))$. Since
\begin{align*}
    \int_{V}\mathbb{E}\|F(t,Y,x)-F(t,\alpha,x)\|^{2}\nu(dx)\leq L\mathbb{E}\|Y-\alpha\|^{2}= L\mathbb{E}\|\beta\|^{2},
    \end{align*}
and $\beta\in SBC_{0}(\mathbb{R},\rho)$, it follows that $F(t,Y(t),x)-F(t,\alpha(t),x)\in PSBC(\mathbb{R}\times V, \rho)$ and
\begin{align*}
\frac{1}{m(r,\rho)}\int_{-r}^{r}\int_{V}\mathbb{E}\|F(t,Y,x)-F(t,\alpha,x)\|^{2}\nu(dx)\rho(t)dt\leq \frac{L}{m(r,\rho)}\int_{-r}^{r}\mathbb{E}\|\beta(t)\|^{2}\rho(t)dt,
  \end{align*}
then we have $F(t,Y(t),x)-F(t,\alpha(t),x)\in PSBC_{0}(\mathbb{R}\times V, \rho)$.

Secondly, we will show that  $\varphi(t,\alpha(t),x)\in PSBC_{0}(\mathbb{R}\times V, \rho)$. Since $\varphi\in PSBC_{0}(\mathbb{R}\times \mathcal{L}^{2}(P,H)\times V, \rho)$, then
$\varphi(t,\alpha(t),x)\in PSBC(\mathbb{R}\times V, \mathcal{L}^{2}(P,H))$. Hence $\varphi(t,\alpha(t),x)$ is Poisson uniformly stochastically continuous on $[-r,r]$. Set $I=\alpha([-r,r])$. We shall show that $I$ is compact in $\mathcal{L}^{2}(P,H))$.
Let $\{\alpha(t_{n})\}$ be an $\mathcal{L}^{2}$ bounded sequence in $I$.
There exists a subsequence $t_{n_{k}}\rightarrow t_{0}$
as $k\rightarrow \infty$.
Since $\alpha$ is $\mathcal{L}^{2}$-continuous, it follows that $\mathbb{E}\|\alpha(t_{n_{k}})-\alpha(t_{0})\|^{2}\rightarrow 0$
as $k\rightarrow \infty$.
Hence, $I$ is a compact set in $\mathcal{L}^{2}(P,H)$. Therefore, one can find finite open balls $O_{k}(k=1,2,\ldots,l)$, with center $\alpha_{k}\in I$ and radius $\delta$ small enough such that $I\subset \bigcup_{k=1}^{l}O_{k}$. By the Poisson uniformly stochastically continuity of $\varphi(t,\alpha(t),x)$, we have for sufficiently small $\delta>0$,
\begin{align}\label{uniform}
\int_{V}\mathbb{E}\|\varphi(t,\alpha(t),x)-F(t,\alpha_{k},x)\|^{2}\nu(dx)<\frac{\varepsilon}{4},
\end{align}
where $\alpha(t)\in U(\alpha_{k}, \delta)$ and $t\in [-r,r]$. \\
The set $B_{k}=\{t\in [-r,r]|\alpha(t)\in O_{k}\}$ is open in $[-r,r]$ and $[-r,r]=\bigcup_{k=1}^{l}B_{k}$. Let
\begin{align*}
E_{1}=B_{1}, \ \ E_{k}=B_{k}\backslash \bigcup_{j=1}^{k-1}B_{k}, \ \ 2\leq k\leq l.
\end{align*}
Then $E_{i}\bigcap E_{j}=\phi$ when $i\neq j, \ 1\leq i,j\leq l$. We have
\begin{align*}
&\{t\in [-r,r]\big|\int_{V}\mathbb{E}\|\varphi(t,\alpha(t),x)\|^{2}\nu(dx)\geq \varepsilon\}\\
&\subset \bigcup_{k=1}^{l}\{t\in E_{k}\big| \int_{U}2\mathbb{E}\|\varphi(t,\alpha(t),x)-\varphi(t,\alpha_{k},x)\|^{2}\nu(dx)+\int_{V}2\mathbb{E}\|\varphi(t,\alpha_{k},x)\|^{2}\nu(dx)\geq \varepsilon\}\\
&\subset \bigcup_{k=1}^{l}\big(\big\{t\in E_{k}\big|\mathbb{E}\|\varphi(t,\alpha(t),x)-\varphi(t,\alpha_{k},x)\|^{2}\nu(dx)\geq \frac{\varepsilon}{4}\}\\
&\cup \{t\in E_{k}\big|\int_{V}\mathbb{E}\|\varphi(t,\alpha_{k},x)\|^{2}\nu(dx)\geq \frac{\varepsilon}{4}\}\big)
\end{align*}
It follows from (\ref{uniform}) that $\{t\in E_{k}\big|\int_{V}\mathbb{E}\|\varphi(t,\alpha(t),x)-F(t,\alpha_{k},x)\|^{2}\nu(dx)\geq\frac{\varepsilon}{4}\}$ are empty for all $1\leq k\leq m$. Therefore
\begin{align*}
\{t\in [-r,r]\big|\int_{V}\mathbb{E}\|\varphi(t,\alpha(t),x)\|^{2}\nu(dx)\geq \varepsilon\}\subset \bigcup_{k=1}^{l}\{t\in E_{k}\big|\int_{V}\mathbb{E}\|\varphi(t,\alpha_{k},x)\|^{2}\nu(dx)\geq \frac{\varepsilon}{4}\}.
\end{align*}
This implies that
\begin{align*}
\frac{1}{m(r,\rho)}\int_{M_{r,\varepsilon}(\varphi(t,\alpha(t),x))}\rho(t)dt\leq \sum_{k=1}^{l}\frac{1}{m(r,\rho)}\int_{M_{r,\frac{\varepsilon}{4}}(\varphi(t,\alpha_{k},x))}\rho(t)dt.
\end{align*}
Note that $\varphi(t,\alpha_{k},x)\in PSBC_{0}(\mathbb{R}\times V,\rho)$. From Lemma 2.32, it follows that
\begin{align*}
\lim_{r\rightarrow +\infty}\frac{1}{m(r,\rho)}\int_{M_{r,\frac{\varepsilon}{4}}(\varphi(t,\alpha_{k},x))}\rho(t)dt=0, \ k=1,2,\ldots,l.
\end{align*}
Thus we have
\begin{align*}
\lim_{r\rightarrow +\infty}\frac{1}{m(r,\rho)}\int_{M_{r,\frac{\varepsilon}{4}}(\varphi(t,\alpha(t),x))}\rho(t)dt=0.
\end{align*}
By Lemma 2.32, $\varphi\in PSBC_{0}(\mathbb{R}\times V, \rho)$. Therefore $J(t,x):=F(t,Y(t),x)\in PSWPAA(\mathbb{R}\times V, \rho)$. \ \ $\Box$\\


\end{document}